\newif \ifNM

\ifNM
\documentclass[a4paper,sn-mathphys]{sn-jnl}
\jyear{2022}%
\theoremstyle{thmstyleone}%
\newtheorem{theorem}{Theorem}
\newtheorem{proposition}[theorem]{Proposition}%
\theoremstyle{thmstyletwo}%
\newtheorem{example}{Example}%
\newtheorem{remark}{Remark}%
\theoremstyle{thmstylethree}%
\newtheorem{definition}{Definition}%
\newtheorem{lemma}{Lemma}%
\newtheorem{prop}{Proposition}%
\else

\documentclass[10pt,a4paper]{article}
\usepackage[a4paper,left=2.5cm,right=2.5cm,top=2.5cm,bottom=2.5cm]{geometry}
\usepackage[normalem]{ulem} 
\usepackage[textsize=footnotesize,color=blue!30!white]{todonotes}
\usepackage{cancel}
\usepackage{pgfplots}

\fi

\usepackage{mathtools}
\usepackage{bm}
\usepackage{amssymb}
\usepackage{amsmath}
\usepackage{enumitem}
\usepackage[normalem]{ulem} 

\ifNM \else

\usepackage{amsthm}
\usepackage[breaklinks,bookmarks=false]{hyperref}

\hypersetup{colorlinks, linkcolor=blue, citecolor=blue,
urlcolor=blue, plainpages=false, pdfwindowui=false, pdfstartview={FitH}, pdftitle={A stable local commuting projector and optimal $hp$ approximation estimates in $H$(curl)}, pdfauthor={Martin Vohral\'ik} }

\fi

\newcommand\cred[1]{%
  \protect\leavevmode
  \begingroup
    #1%
  \endgroup
}

\numberwithin{equation}{section}

\newif \ifbbm

\newif \iffract

\newcommand{\be}{\begin{equation}}
\newcommand{\ee}{\end{equation}}
\newcommand{\bea}{\begin{eqnarray}}
\newcommand{\eea}{\end{eqnarray}}
\newcommand{\bean}{\begin{eqnarray*}}
\newcommand{\eean}{\end{eqnarray*}}
\def\ba#1\ea{\begin{align}#1\end{align}}
\def\ban#1\ean{\begin{align*}#1\end{align*}}
\def\bat#1\eat{\begin{alignat}#1\end{alignat}}
\def\batn#1\eatn{\begin{alignat*}#1\end{alignat*}}
\def\bs#1\es{\begin{split}#1\end{split}}
\newcommand{\bse}{\begin{subequations}}
\newcommand{\ese}{\end{subequations}}
\newcommand{\bt}{\begin{theorem}}
\newcommand{\et}{\end{theorem}}
\newcommand{\bpr}{\begin{proposition}}
\newcommand{\epr}{\end{proposition}}
\newcommand{\bcj}{\begin{conjecture}}
\newcommand{\ecj}{\end{conjecture}}
\newcommand{\bl}{\begin{lemma}}
\newcommand{\el}{\end{lemma}}
\newcommand{\bc}{\begin{corollary}}
\newcommand{\ec}{\end{corollary}}
\newcommand{\bp}{\begin{proof}}
\newcommand{\ep}{\end{proof}}
\newcommand{\bd}{\begin{definition}}
\newcommand{\ed}{\end{definition}}
\newcommand{\brem}{\begin{remark}}
\newcommand{\erem}{\end{remark}}
\newcommand{\bas}{\begin{assumption}}
\newcommand{\eas}{\end{assumption}}
\newcommand{\bex}{\begin{example}}
\newcommand{\eex}{\end{example}}
\newcommand{\bqo}{\begin{quote}}
\newcommand{\eqo}{\end{quote}}
\newcommand{\bdc}{\begin{description}}
\newcommand{\edc}{\end{description}}
\newcommand{\bi}{\begin{itemize}}
\newcommand{\ei}{\end{itemize}}
\newcommand{\ben}{\begin{enumerate}}
\newcommand{\een}{\end{enumerate}}

\ifNM
\newtheorem{assumption}[theorem]{Assumption}

\else

\newtheorem{lemma}{Lemma}[section]
\newtheorem{corollary}[lemma]{Corollary}
\newtheorem{prop}[lemma]{Proposition}
\newtheorem{assumption}[lemma]{Assumption}
\newtheorem{theorem}[lemma]{Theorem}
\newtheorem{definition}[lemma]{Definition}
\newtheorem{remark}[lemma]{Remark}
\newtheorem{example}[lemma]{Example}

\newtheorem{conjecture}[lemma]{Conjecture}

\fi

\newcommand\Om{\Omega}
\newcommand\om{\omega}
\newcommand\oma{{\omega_\ver}} 
\newcommand\omb{{\omega_\vertt}}

\newcommand\omF{{\omega_\sd}} 
\newcommand\GD{{\Gamma_{\mathrm D}}} 
\newcommand\GN{{\Gamma_{\mathrm N}}} 
\newcommand\gD{{\gamma_{\mathrm D}}} 

\newcommand\Gr{\nabla} 
\newcommand\Dv{\nabla {\cdot}} 
\newcommand\Crl{\nabla {\times}} 

\newcommand\gr{\mathrm{grad}}
\newcommand\dv{\mathrm{div}}

\newcommand\crl{\mathrm{curl}}
\newcommand\scp{{\cdot}} 
\newcommand\vp{{\times}} 
\newcommand{\jump}[1]{[\![#1]\!]}

\DeclarePairedDelimiter\norm{\|}{\|} 

\ifbbm

\else

\fi

\newcommand\HO{H^1}
\newcommand\Ho{H^1(\Om)}

\newcommand\tHo{{\bm H}^1(\Om)}
\newcommand\Hoi[1]{H^1(#1)}
\newcommand\Hti[1]{H^2(#1)}

\newcommand\tHoi[1]{{\bm H}^1(#1)}
\newcommand\tHoDi[1]{{\bm H}^1_{0,\mathrm{D}}(#1)}

\newcommand\HoD{H^1_{0,\mathrm{D}}(\Om)}
\newcommand\HoN{H^1_{0,\mathrm{N}}(\Om)}

\newcommand\Hooi[1]{H^1_0(#1)}

\newcommand\Lt{L^2(\Om)}
\newcommand\Lts{L^2_*(\Om)}
\newcommand\LT{L^2}
\newcommand\tLt{\bm{L}^2(\Om)}
\newcommand\tLti[1]{\bm{L}^2(#1)}
\newcommand\tLT{\bm{L}^2}
\newcommand\Lti[1]{L^2(#1)}

\newcommand\HDV{\bm{H}(\dv)}
\newcommand\Hdv{\bm{H}(\dv,\Om)}
\newcommand\HdvN{\bm{H}_{0,\mathrm{N}}(\dv,\Om)}

\newcommand\Hdvi[1]{\bm{H}(\dv,#1)}

\newcommand\Hdva{\bm{H}_{0}(\dv, \oma)}

\newcommand\HC{\bm{H}(\crl)}
\newcommand\Hc{\bm{H}(\crl, \Om)}
\newcommand\Hci[1]{\bm{H}(\crl, #1)}
\newcommand\Hciz[1]{\bm{H}_0(\crl, #1)}

\newcommand\HcD{\bm{H}_{0,\mathrm{D}}(\crl, \Om)}
\newcommand\HcDi[1]{\bm{H}_{0,\mathrm{D}}(\crl, #1)}
\newcommand\HcN{\bm{H}_{0,\mathrm{N}}(\crl, \Om)}
\newcommand\Hca{\bm{H}_{0}(\crl, \oma)}

\newcommand\Hcdual[1]{\bm{H}^\dagger(\crl, #1)}

\newcommand\RT{\bm{\mathcal{R\hspace{-0.1em}T}}\hspace{-0.25em}}
\newcommand\RTproj[1]{\bm{I}_{\bm{\mathcal{R\hspace{-0.1em}T}}}^{#1}}
\newcommand\RToproj[1]{\bm{\Pi}_{\bm{\mathcal{R\hspace{-0.1em}T}}}^{#1}}
\newcommand\ND{\bm{\mathcal{N}}\hspace{-0.2em}}
\newcommand\NDproj[1]{\bm{I}_{\bm{\mathcal{N}}}^{#1}}

\newcommand\Phc{\bm{P}_h^{p,\crl}}
\newcommand\Phd{\bm{P}_h^{p,\dv}}

\ifbbm

\else

\fi

\newcommand\ie{i.e.}
\newcommand\cf{cf.}
\newcommand\eg{e.g.}
\newcommand\eal{{\em et al.}}

\newcommand\eq{:=}

\newcommand\ls{\lesssim}

\newcommand\nn{\nonumber}
\newcommand\pt{\partial}
\newcommand{\<}{\langle}
\renewcommand{\>}{\rangle}

\newcommand{\suma}{\sum_{\ver \in \Vh}}

\newcommand\nv{\bm 0}
\newcommand\reff[2]{\stackrel{\eqref{#1}}{#2}}

\newcommand\refff[3]{\stackrel{{\substack{\eqref{#1}\\ \eqref{#2}}}}{#3}}
\newcommand\refvs[3]{\stackrel{\mathmakebox[\widthof{#3}]{#1}}{#2}}

\ifbbm

\else

\fi

\newcommand{\elm}{{K}} 
\newcommand{\elmtt}{{K'}}
\newcommand{\sd}{{F}} 

\newcommand{\ver}{{\bm{a}}} 
\newcommand{\vertt}{{\bm{a}'}}

\newcommand\FK{\mathcal{F}_\elm}

\newcommand\Faint{\mathcal{F}_{\ver}^\mathrm{int}}

\newcommand\EK{\mathcal{E}_\elm}

\newcommand\Th{{\mathcal{T}_h}} 
\newcommand\TK{\mathcal{T}_\elm} 
\newcommand\Ta{{\mathcal{T}_{\ver}}} 

\newcommand\Vh{\mathcal{V}_h} 

\newcommand\VK{\mathcal{V}_\elm}
\newcommand\VKtt{\mathcal{V}_\elmtt}

\newcommand\psia{\psi^\ver}
\newcommand\psiF{\psi^\sd}

\newcommand\psib{\psi^\vertt}

\newcommand\frh{\bsig_h} 

\newcommand\btha{\bm{\theta}_h^\ver}
\newcommand\bthb{\bm{\theta}_h^\vertt}
\newcommand\bdla{\bdl_h^\ver}
\newcommand\bha{\bm{h}_h^\ver}

\ifbbm

\else

\fi

\newcommand\vf{\varphi}

\newcommand\tx{{\bm{x}}}

\newcommand\td{\bm{d}}
\newcommand\te{\bm{e}}
\newcommand\tf{\bm{f}}

\newcommand\tm{\bm{m}}
\newcommand\tn{\bm{n}}
\newcommand\tp{\bm{p}}

\newcommand\tr{\bm{r}}

\newcommand\tv{\bm{v}}
\newcommand\ttw{\bm{w}}

\newcommand\ttau{\bm{\tau}}

\newcommand\tL{\bm{L}}

\newcommand{\bPi}{\bm{\Pi}}
\newcommand{\bsig}{\bm{\sigma}}

\newcommand{\bdl}{\bm{\delta}}
\newcommand\bvf{\bm{\varphi}}

\newcommand\bpsi{\bm{\psi}}
\newcommand\bchi{\bm{\chi}}
\newcommand\btau{\bm{\tau}}
\newcommand\biot{\bm{\iota}}
\newcommand\bh{\bm{h}}
\newcommand\bzeta{\bm{\zeta}}

\newcommand\RR{{\mathbb R}}

\ifbbm

\newcommand\PP{{\mathbb P}}

\else

\newcommand\PP{{\mathcal P}}

\fi

\iffract

\newcommand\ft{{\frac 1 2}}

\else

\newcommand\ft{{1/2}}

\fi

\newcommand\CPFi[1]{C_{\mathrm{PF},#1}}
\newcommand\Clift[1]{C_{\mathrm{L},#1}}

\newcommand\tHsa{{\bm H}^1_*(\oma)}
\newcommand\Hsa{H^1_*(\oma)}

\def \abs{We design an operator from the infinite-dimensional Sobolev space $\HC$ to its finite-dimensional subspace formed by the N\'ed\'elec piecewise polynomials on a tetrahedral mesh that has the following properties: 1) it is defined over the entire $\HC$, including boundary conditions imposed on a part of the boundary; 2) it is defined locally in a neighborhood of each mesh element; 3) it is based on simple piecewise polynomial projections; 4) it is stable in the $\tL^2$-norm, up to data oscillation; 5) it has optimal (local-best) approximation properties; 6) it satisfies the commuting property with its sibling operator on $\HDV$; 7) it is a projector, \ie, it leaves intact objects that are already in the N\'ed\'elec piecewise polynomial space. This operator can be used in various parts of numerical analysis related to the $\HC$ space. We in particular employ it here to establish the two following results: i) equivalence of global-best, tangential-trace- and curl-constrained, and local-best, unconstrained approximations in $\HC$ including data oscillation terms; and ii) fully $h$- and $p$- (mesh-size- and polynomial-degree-) optimal approximation bounds valid under the minimal Sobolev regularity only requested elementwise. As a result of independent interest, we also prove a $p$-robust equivalence of curl-constrained and unconstrained best-approximations on a single tetrahedron in the $\HC$-setting, including $hp$ data oscillation terms.}

\def\KW{Sobolev space $\HC$, N\'ed\'elec finite element space, stable local commuting projector, constrained--unconstrained equivalence, local-best approximation, $hp$ approximation, minimal local Sobolev regularity}

\ifNM

\title[A stable local commuting projector and $hp$ estimates in $\HC$]{A stable local commuting projector and optimal $hp$ approximation estimates in $\HC$}

\author[1]{\fnm{Th\'eophile} \sur{Chaumont-Frelet}}\email{theophile.chaumont@inria.fr}

\author*[2,3]{\fnm{Martin} \sur{Vohral\'ik}}\email{martin.vohralik@inria.fr}

\affil[1]{\orgname{University C\^ote d'Azur, Inria, CNRS, LJAD}, \orgaddress{\street{2004 Route des Lucioles}, \city{Valbonne}, \postcode{06902}, \country{France}}}

\affil[2]{\orgdiv{SERENA}, \orgname{Inria}, \orgaddress{\street{2 rue Simone Iff}, \city{Paris}, \postcode{75589}, \country{France}}}

\affil[3]{\orgdiv{CERMICS}, \orgname{Ecole des Ponts}, \orgaddress{\street{6 et 8 avenue Blaise Pascal}, \city{Marne-la-Vall\'ee}, \postcode{77455}, \country{France}}}

\abstract{\abs}


\else

\title{A stable local commuting projector and optimal $hp$ approximation estimates in $\HC$\thanks{This project has received funding from the European Research Council (ERC) under the European Union's Horizon 2020 research and innovation program (grant agreement No 647134 GATIPOR).}}

\author{Th\'eophile Chaumont-Frelet\footnotemark[2] \and Martin Vohral\'ik\footnotemark[3]\,\,\footnotemark[4]}

\fi

\begin{document}
\maketitle

\ifNM\else

\renewcommand{\thefootnote}{\fnsymbol{footnote}}

\footnotetext[2]{University C\^ote d'Azur, Inria, CNRS, LJAD, 2004 Route des Lucioles, 06902 Valbonne, France
(\href{mailto:theophile.chaumont@inria.fr}{\texttt{theophile.chaumont@inria.fr}}).}
\footnotetext[3]{Inria, 2 rue Simone Iff, 75589 Paris, France
(\href{mailto:martin.vohralik@inria.fr}{\texttt{martin.vohralik@inria.fr}}).}
\footnotetext[4]{CERMICS, Ecole des Ponts, 77455 Marne-la-Vall\'ee, France}

\renewcommand{\thefootnote}{\arabic{footnote}}

\begin{abstract}
\abs
\end{abstract}

\fi

\ifNM
\else
\fi

\ifNM
\else


\fi

\section{Introduction} \label{sec_intr}

Let $\Om \subset \RR^3$ be a Lipschitz polyhedral domain (open, bounded, and connected set) and $\GN$ a Lipschitz polygonal relatively open subset of its boundary $\pt \Om$ (details on setting and notation are given in Section~\ref{sec_not} below). A central concept in numerical analysis of
partial differential equations including the $\gr$, $\crl$, and $\dv$ operators, connected with
the Sobolev spaces $\HoN$, $\HcN$, and $\HdvN$, is the following commuting de Rham complex:
\ifNM
\be \label{eq_ex_seq}\arraycolsep=2.2pt
    \begin{array}{ccccccc}
      \HoN & \xrightarrow{\,\Gr\,} & \HcN & \xrightarrow{\,\Crl\,} & \HdvN & \xrightarrow{\,\Dv\,} & \Lts\\
      \Big \downarrow P_h^{\cred{p+1},\gr} &  & \Big \downarrow \Phc &  & \Big \downarrow \Phd & & \Big \downarrow \Pi_h^p  \\
      \PP_{\cred{p+1}}(\Th)\cap & \xrightarrow{\,\Gr\,} & \ND_p(\Th) \cap & \xrightarrow{\,\Crl\,} & \RT_p(\Th) \cap & \xrightarrow{\,\Dv\,} & \PP_p(\Th) \cap \\
      \HoN & & \HcN & & \HdvN & & \Lts
    \end{array}
\ee
\else
\be \label{eq_ex_seq}\arraycolsep=2.2pt
    \begin{array}{ccccccc}
      \HoN & \xrightarrow{\,\Gr\,} & \HcN & \xrightarrow{\,\Crl\,} & \HdvN & \xrightarrow{\,\Dv\,} & \Lts\\
      \Big \downarrow P_h^{\cred{p+1},\gr} &  & \Big \downarrow \Phc &  & \Big \downarrow \Phd & & \Big \downarrow \Pi_h^p  \\
      \PP_p(\Th)\cap \HoN & \xrightarrow{\,\Gr\,} & \ND_p(\Th) \cap \HcN & \xrightarrow{\,\Crl\,} & \RT_p(\Th) \cap \HdvN & \xrightarrow{\,\Dv\,} & \PP_p(\Th) \cap \Lts.\\
    \end{array}
\ee
\fi

Assuming for simplicity in the introduction that $\Omega$
\cred{is homotopic to a ball (contractible) and either $\GD = \pt \Om$ or $\GN = \pt \Om$
so that the cohomology spaces are trivial,} the first line of~\eqref{eq_ex_seq} is the
well-known exact sequence on the continuous, infinite-dimensional, level, see, \eg,
Arnold~\eal\ \cite{Arn_Falk_Winth_FEC_06} and the references therein.
It in particular states that 1) each function from the $\HcN$ space whose weak curl
vanishes is a weak gradient of a function from $\HoN$; 2) each function from the $\HdvN$
space whose weak divergence vanishes is a weak curl of a function from $\HcN$; 3) each function
from the $\Lts$ space is a weak divergence of a function from $\HdvN$.
Similarly, the second line is the counterpart of the first one on the discrete,
finite-dimensional, piecewise polynomial, level, see \eg, Boffi~\eal\ \cite{Bof_Brez_For_MFEs_13}
and the references therein. The passage between the first and the second line is then the key
interest in this contribution, where the three operators $P_h^{\cred{p+1},\gr}$, $\Phc$, and
$\Phd$\cred{, $p \geq 0$,} should:

\ben

\item \label{pr_ent_sp} be defined over the {\em entire infinite-dimensional spaces} $\HoN$, $\HcN$, and $\HdvN$;

\item be defined {\em locally}, in a neighborhood of mesh elements at most;

\item be based on {\em simple} piecewise polynomial projections;

\item be {\em stable} in $\tLt$ for $\HcN$ and $\HdvN$ and in $\tLt$ of the weak gradient for $\HoN$, up to data oscillation;

\item \label{pr_appr} have {\em optimal approximation properties}, i.e., that of {\em local-best} unconstrained $\tLT$-orthogonal projectors;

\item satisfy the {\em commuting properties} expressed by the arrows in~\eqref{eq_ex_seq};

\item \label{pr_proj}  be {\em projectors}, \ie, leave intact objects that are already in the piecewise polynomial spaces.

\een

There is an immense literature devoted to~\eqref{eq_ex_seq}. A first consideration for the operators $P_h^{\cred{p+1},\gr}$, $\Phc$, and $\Phd$ is given by the canonical projectors, see Ciarlet~\cite{Ciar_78}, N\'ed\'elec~\cite{Ned_mix_R_3_80}, and Raviart and Thomas~\cite{Ra_Tho_MFE_77}, respectively. These satisfy many of the properties above, but, unfortunately, not property~\ref{pr_ent_sp}, since their action is not defined on all objects from the entire infinite-dimensional spaces $\HoN$, $\HcN$, and $\HdvN$. The commuting diagram~\eqref{eq_ex_seq} has been addressed at the abstract level of the finite element exterior calculus in, \eg, Christiansen and Winther~\cite{Christ_Wint_sm_proj_08}, still leading to the loss of some of the desirable properties, namely the locality. Simultaneous definition on the entire Sobolev spaces, locality, commutativity, and the projection property have been achieved in Falk and Winther~\cite{Falk_Winth_loc_coch_14}, though the stability in the $\tLt$ norms (up to data oscillation) and the local-best unconstrained approximation properties have not been addressed; the stability in the $\tLt$ norms been recently achieved in Arnold and Guzm\'{a}n~\cite{Arn_Guz_loc_stab_L2_com_proj_21}. Two different sets of projectors, satisfying together (but not individually) all properties~\ref{pr_ent_sp}--\ref{pr_proj}, were then designed in Ern and Guermond~\cite{Ern_Guer_mol_de_Rham_16, Ern_Guer_quas_int_best_appr_17}. Finally
an operator $\Phd$ satisfying the integrality of the requested properties
\cred{(all~\ref{pr_ent_sp}--\ref{pr_proj} above)} has been recently devised in Ern~\eal\ \cite[Section~3.1]{Ern_Gud_Sme_Voh_loc_glob_div_22}.

The first goal of the present contribution is to design an operator $\Phc$ satisfying the {\em integrality of the requested properties}~\ref{pr_ent_sp}--\ref{pr_proj}. Definition~\ref{def_Phc} is designed to this purpose, relying on (a slight modification of) $\Phd$ from~\cite{Ern_Gud_Sme_Voh_loc_glob_div_22}, see Definition~\ref{def_Phd}, and using similar building principles as in~\cite{Ern_Gud_Sme_Voh_loc_glob_div_22}. The main result here is Theorem~\ref{thm_Phc}.
The central technical tool allowing to achieve the commuting property is related to equilibration in $\Hc$. A first contribution in this direction is that of Braess and Sch{\"o}berl~\cite{Braess_Scho_a_post_edge_08}. Recent extensions to higher polynomial degrees are developed in Gedicke~\eal\ \cite{Ged_Gee_Per_a_post_Maxw_20, Ged_Gee_Per_Sch_post_Maxw_21} as well as in~\cite{Chaum_Voh_Maxwell_equil_23}, which we use here.
\cred{All~\ref{pr_ent_sp}--\ref{pr_proj} come handy namely in a priori error estimates (especially property~\ref{pr_appr}) and in approximation classes in proofs of convergence and optimality based on posteriori error estimates. All~\ref{pr_ent_sp}--\ref{pr_proj} also turn useful in localized orthogonal decomposition techniques for multiscale problems, \cf\ \cite{Gal_Henn_Verf_num_hom_H_curl_18} and the references therein.}

Our contribution stands apart from the existing literature namely in the satisfaction of property~\ref{pr_appr}. This leads to the result of equivalence of global-best (tangential-trace- and curl-constrained) and local-best (unconstrained) approximations in $\Hc$, see Theorem~\ref{thm_loc_glob}. This result, not taking into account data oscillation, has been recently established in~\cite{Chaum_Voh_loc_glob_curl_21}, building on the seminal contribution by Veeser~\cite{Veeser_approx_grads_16} in the $\Ho$-setting and on~\cite{Ern_Gud_Sme_Voh_loc_glob_div_22} in the $\Hdv$-setting, see also the references therein. Here, we present a direct proof. We take into account data oscillation, which actually turns out quite demanding.
\cred{A related result giving rise to an operator with local-best approximation properties in the general framework of finite element spaces of differential forms has recently been obtained in Gawlik~\eal\ \cite{Gaw_Holst_Licht_loc_approx_FEEC_21}; it, however, does not address the curl constraint.}

Yet a separate, and involved, question in numerical analysis is that of deriving
$hp$-approximation estimates. This has been addressed in the $\Hdv$- and $\Hc$-settings
in particular in Suri~\cite{Sur_stab_cvg_MFE_p_90}, Monk~\cite{Monk_p_hp_Nedelec_94},
Demkowicz and Buffa~\cite{Demk_Buf_q_opt_proj_int_05}, and Demkowicz~\cite{Demk_pol_ex_seq_int_06},
see also the references therein. These references feature a slight suboptimality in the polynomial
degree $p$ on tetrahedral meshes (presence of a logarithmic factor), which has been removed in Bespalov and Heuer~\cite{Besp_Heue_H_curl_hp_2D_09, Besp_Heue_H_div_hp_2D_11}
and recently in Melenk and Rojik~\cite{Mel_Roj_com_p_interp_20}. Unfortunately, none of these references allows for minimal Sobolev regularity.
The result in~\cite[Theorem~3.6]{Ern_Gud_Sme_Voh_loc_glob_div_22} is equally fully
$h$- and $p$-optimal, and this, moreover, under the minimal Sobolev regularity, only requested elementwise.
Deriving such estimates in the $\Hc$-setting is the last goal of the present contribution.
Theorem~\ref{thm_hp} in particular presents a fully $h$- and $p$-optimal approximation bound
valid under the minimal Sobolev regularity that is only requested separately on each mesh
element (\cred{the treatment of the $hp$ data oscillation
is restricted to convex patches or presents a slight suboptimality}).  The key ingredients are here again linked to (polynomial-degree-robust) flux
equilibration in $\Hc$ of~\cite{Chaum_Voh_Maxwell_equil_23}, with the cornerstones being
the results in a single tetrahedron: on the right inverses \cred{in the volume} by
Costabel and McIntosh~\cite{Cost_McInt_Bog_Poinc_10} and on the polynomial extension operators
\cred{from the boundary} by Demkowicz~\eal\
\cite{Demk_Gop_Sch_ext_I_09, Demk_Gop_Sch_ext_II_09, Demk_Gop_Sch_ext_III_12}.

This contribution is organized as follows: Section~\ref{sec_not} fixes the setting and notation.
The above-described main results are collected in Section~\ref{sec_main_res}. The well-posedness
of the central definition of our stable local commuting projector is verified in
Section~\ref{sec_Phc_WP}, and the proofs of the three principal theorems are then presented
respectively in Sections~\ref{sec_Phc_proof}--\ref{sec_hp_proof}. A result of independent
interest, stipulating a polynomial-degree-robust equivalence of constrained and unconstrained
best-approximation on a single tetrahedron in the $\Hc$-setting, including $hp$ data oscillation
terms, is presented in Appendix~\ref{sec_constr}. \cred{Finally, a technical result on broken
regular decomposition in a patch is presented in Appendix~\ref{sec_broken_reg}.}

\section{Setting and notation} \label{sec_not}

Let $\om, \Om \subset \RR^3$ be open, Lipschitz polyhedral domains; $\Om$ will be
used to denote the computational domain, while we reserve the
notation $\om \subseteq \Om$ for its subdomains. \cred{We do not make any assumptions
on the topology of $\Om$; $\Om$ is not necessarily homotopic to a ball (contractible)
and nontrivial cohomology spaces induced by $\Om$ are allowed. The subdomains $\om$
(patches of elements below) will later, in turn, be supposed homotopic to a ball (contractible).}
We will use the notation $a \ls b$ when there exists a positive constant $C$ such that
$a \leq C b$; we will always specify the dependencies of $C$.

\subsection{Sobolev spaces \texorpdfstring{$\HO$, $\HC$, and $\HDV$}{H1, H(curl), and H(div)}} \label{sec_bas_cont_sp}

We let $\Lti{\om}$ be the space of scalar-valued square-integrable functions defined on $\om$;
we use the notation $\tLti{\om} \eq [\Lti{\om}]^3$ for vector-valued functions with each
component in $\Lti{\om}$. We denote by $\norm{{\cdot}}_\om$ the $\Lti{\om}$ or $\tLti{\om}$
norm and by $({\cdot},{\cdot})_\om$ the corresponding scalar product; we drop the index when
$\om = \Om$. We will extensively work with the following three Sobolev spaces: 1) $\Hoi{\om}$,
the space of scalar-valued $\Lti{\om}$ functions with weak gradients in $\tLti{\om}$,
$\Hoi{\om} \eq \{v \in \Lti{\om}; \, \Gr v \in \tLti{\om}\}$; 2) $\Hci{\om}$,
the space of vector-valued $\tLti{\om}$ functions with weak curls in $\tLti{\om}$,
$\Hci{\om} \eq \{\tv \in \tLti{\om}; \, \Crl \tv \in \tLti{\om}\}$; and 3) $\Hdvi{\om}$,
the space of vector-valued $\tLti{\om}$ functions with weak divergences in $\Lti{\om}$,
$\Hdvi{\om} \eq \{\tv \in \tLti{\om}; \, \Dv \tv \in \Lti{\om}\}$. We refer the reader to
Adams~\cite{Adams_75} and Girault and Raviart~\cite{Gir_Rav_NS_86} for an in-depth description
of these spaces. Moreover, component-wise $\Hoi{\om}$ functions will be denoted by
$\tHoi{\om} \eq \{\tv \in \tLti{\om}; \, \tv_i \in \Hoi{\om}, \, i=1, \ldots, 3\}$.
We will employ the notation $\<{\cdot},{\cdot}\>_S$ for the integral product on boundary (sub)sets $S \subset \pt \om$.

\subsection{Tetrahedral mesh, patches of elements, and the hat functions} \label{sec_mesh}

Let $\Th$ be a simplicial mesh of the domain $\Om$, \ie,
$\cup_{\elm \in \Th} \elm = \overline \Om$, where any element $\elm \in \Th$ is a
closed tetrahedron with nonzero measure, and where the intersection of two
different tetrahedra is either empty or their common vertex, edge, or face. The shape-regularity
parameter of the mesh $\Th$ is the positive real number
$\kappa_\Th \eq \max_{\elm \in \Th} h_\elm / \rho_\elm$, where $h_\elm$ is the diameter of the
tetrahedron $\elm$ and $\rho_\elm$ is the diameter of the largest ball contained in $\elm$.
These assumptions are standard, and allow for strongly graded meshes with local refinements, though not for anisotropic elements.

We denote the set of vertices of the mesh $\Th$ by $\Vh$; it is composed of interior vertices
lying in $\Om$ and of vertices lying on the boundary $\pt \Om$.
For an element $\elm \in \Th$, $\FK$ denotes the set of its faces and $\VK$ the set of its vertices.
Conversely, for a vertex $\ver \in \Vh$, $\Ta$ denotes the patch of the elements of $\Th$ that share $\ver$, and $\oma$ is the corresponding open subdomain with diameter $h_\oma$.
\cred{We suppose these vertex patch subdomains $\oma$ to be homotopic to a ball (contractible).}

A particular role below will be played by the continuous, piecewise affine ``hat'' function $\psia$ which takes value $1$ at the vertex $\ver$ and zero at the other vertices.
We note that $\oma$ corresponds to the support of $\psia$ and that the functions $\psia$ form the partition of unity
\be \label{eq_PU}
    \suma \psia = 1.
\ee
%
By $\jump{v}$, we denote the jump of the function $v$ on a face $\sd$, \ie,
the difference of the traces of $v$ from the two elements sharing $\sd$
along an arbitrary but fixed normal.

\subsection{Sobolev spaces with partially vanishing traces on \texorpdfstring{$\Om$}{Omega} and \texorpdfstring{$\oma$}{omega}} \label{sec_cont_sp_BC}

Let $\GD$, $\GN$ be two disjoint, relatively open, and possibly empty subsets of the
computational domain boundary $\pt \Om$ such that $\pt \Om = \overline \GD \cup \overline \GN$.
We also require that $\GD$ and $\GN$ have polygonal Lipschitz boundaries, and
we assume that each boundary face of the mesh $\Th$ lies entirely either
in $\overline \GD$ or in $\overline \GN$.
\cred{Notice that the assumption that $\GD$ and $\GN$ have Lipschitz boundaries
excludes ``checkerboard'' patterns of mixed boundary conditions. In particular,
\cred{the cohomology spaces with the boundary conditions on $\oma$ introduced below are trivial}.}

\cred{We denote by} $\Lts$ the subspace of $\Lt$ \cred{formed by} functions of mean value $0$
if $\GN = \pt \Om$ \cred{and set $\Lts = \Lt$ otherwise.} $\HoD$ is the subspace of $\Ho$
formed by functions vanishing on $\GD$ in the sense of traces,
$\HoD \eq \{v \in \Ho; \, v=0$ on $\GD\}$.
Let $\tn_\Om$ be the unit normal vector on $\pt \Om$, outward to $\Om$.
Then $\HcN$, $\HcD$ are the subspaces of $\Hc$ formed by
functions with vanishing tangential traces respectively on $\GN$ or $\GD$,
\cred{\bse \label{eq_HcN_HcD}
\ba
    \HcN &\eq \{\tv \in \Hc; \, \tv \vp \tn_{\Om}=0 \text{ on } \GN\}, \label{eq_HcN} \\
    \HcD &\eq \{\tv \in \Hc; \, \tv \vp \tn_{\Om}=0 \text{ on }\GD\}, \label{eq_HcD}
\ea \ese
where
\be \label{eq_tang_trace} \bs
    \tv \vp \tn_{\Om} = 0 \text{ on } \GN \Longleftrightarrow {} & (\Crl \tv, \bvf) - (\tv, \Crl \bvf) = 0 \\
    {} & \qquad \forall \bvf \in \tHo \text{ such that } \bvf \vp \tn_{\Om} = \nv \text{ on } \GD
\es \ee
and symmetrically for $\GD$.}
Finally, $\HdvN$ is the subspace of $\Hdv$ formed by functions with vanishing normal
trace on $\GN$,
\be \label{eq_HdvN}
    \cred{\HdvN \eq \{\tv \in \Hdv; \, \tv \scp \tn_{\Om}=0 \text{ on } \GN\}},
\ee
where
\be \label{eq_norm_trace}
    \cred{\tv \scp \tn_{\Om} = 0 \text{ on } \GN \Longleftrightarrow (\tv, \Gr \vf) + (\Dv \tv, \vf) = 0 \qquad \forall \vf \in \HoD}.
\ee
Fernandes and Gilardi~\cite{Fer_Gil_Maxw_BC_97} present a
rigorous characterization of tangential (resp. normal) traces of $\Hc$ (resp. $\Hdv$) on
a part of the boundary $\pt \Om$.

We will also need local spaces on the patch subdomains $\oma$. Let first $\ver \in \Vh$ be
an interior vertex. Then we \cred{set
\bse \label{eq_loc_spaces_int} \ba
    \Hsa & \eq \{v \in \Hoi{\oma}; \, (v,1)_\oma = 0\}, \label{eq_Hsa} \\
    \Hca & \eq \{\tv \in \Hci{\oma}; \, \tv \vp \tn_{\oma}=0 \text{ on } \pt \oma\}, \label{eq_Hca}\\
    \Hcdual{\oma} & \eq \Hci{\oma}, \label{eq_Hca_dual}\\
    \Hdva & \eq \{\tv \in \Hdvi{\oma}; \, \tv \scp \tn_{\oma}=0 \text{ on } \pt \oma\}, \label{eq_Hdva}
\ea \ese
where the tangential trace in~\eqref{eq_Hca} is understood by duality as above in~\eqref{eq_tang_trace}, whereas the normal trace in~\eqref{eq_Hdva} is understood by duality as above in~\eqref{eq_norm_trace}. By definition, $\Hsa$ from~\eqref{eq_Hsa} is the subspace of those $\Hoi{\oma}$ functions whose mean value vanishes.}

The situation is more subtle for boundary vertices. As a first possibility,
if $\ver \in \GN$ (\ie, $\ver \in \Vh$ is a boundary vertex such that all the boundary faces sharing the vertex $\ver$ lie in $\overline{\GN}$), then the spaces $\Hsa$, $\Hca$, $\Hcdual{\oma}$, and $\Hdva$ are defined as above \cred{in~\eqref{eq_loc_spaces_int}}.
Secondly, when $\ver \in \overline \GD$, then at least one of the faces
sharing the vertex $\ver$ lies in $\overline \GD$, and we denote by $\gD$ the subset of $\GD$
corresponding to all such faces. \cred{Notice that due to our assumptions on $\GD$
and $\GN$ to have polygonal Lipschitz boundaries, $\gD$ is always simply connected.}
In this situation, we \cred{let
\bse \label{eq_loc_spaces_Dir} \ba
    \Hsa & \eq \{v \in \Hoi{\oma}; \, v = 0 \text{ on } \gD\}, \label{eq_Hsa_Dir} \\
    \Hca & \eq \{\tv \in \Hci{\oma}; \, \tv \vp \tn_{\oma}=0 \text{ on } \pt \oma \setminus \gD\}, \label{eq_Hca_Dir}\\
    \Hcdual{\oma} & \eq \{\tv \in \Hci{\oma}; \, \tv \vp \tn_{\oma}=0 \text{ on } \gD\}, \label{eq_Hca_dual_Dir}\\
    \Hdva & \eq \{\tv \in \Hdvi{\oma}; \, \tv \scp \tn_{\oma}=0 \text{ on } \pt \oma \setminus \gD\}. \label{eq_Hdva_Dir}
\ea \ese}
In all cases, component-wise $\Hsa$ functions are denoted by $\tHsa$.

\subsection{Piecewise polynomial spaces}
\label{sec_disc_sp}

Let $q \geq 0$ be an integer. For a single tetrahedron $\elm \in \Th$, we denote
by $\PP_q(\elm)$ the space of scalar-valued polynomials on $\elm$ of total degree at most $q$,
and by $[\PP_q(\elm)]^3$ the space of vector-valued polynomials on $\elm$ with each
component in $\PP_q(\elm)$. The N\'ed\'elec~\cite{Bof_Brez_For_MFEs_13,Ned_mix_R_3_80} space of
degree $q$ on $\elm$ is then given by
\be \label{eq_N_K}
    \ND_q(\elm) \eq [\PP_q(\elm)]^3 + \tx \vp [\PP_q(\elm)]^3.
\ee
Similarly, the Raviart--Thomas~\cite{Bof_Brez_For_MFEs_13,Ra_Tho_MFE_77} space of degree $q$
on $\elm$ is given by
\be \label{eq_RT_K}
    \RT_q(\elm) \eq [\PP_q(\elm)]^3 + \PP_q(\elm) \tx.
\ee
We note that~\eqref{eq_N_K} and~\eqref{eq_RT_K} are equivalent to the writing with a direct
sum and only homogeneous polynomials in the second terms.
The second term in~\eqref{eq_N_K} is also equivalently given by homogeneous $(q+1)$-degree polynomials $\tv_h$ such that $\tx \scp \tv_h(\tx) = 0$ for all $\tx \in \elm$.

We will below extensively use the broken, piecewise polynomial spaces formed from the above
element spaces
\begin{align*}
    \PP_q(\Th) & \eq \{  v_h \in  \Lt; \,  v_h\vert_\elm \in \PP_q(\elm)\quad \forall \elm \in \Th\},\\
    \ND_q(\Th) & \eq \{\tv_h \in \tLt; \, \tv_h\vert_\elm \in \ND_q(\elm)\quad \forall \elm \in \Th\},\\
    \RT_q(\Th) & \eq \{\tv_h \in \tLt; \, \tv_h\vert_\elm \in \RT_q(\elm)\quad \forall \elm \in \Th\}.
\end{align*}
To form the usual finite-dimensional Sobolev subspaces, we will write $\PP_q (\Th) \cap \Ho$
(for $q \geq 1$), $\ND_q(\Th) \cap \Hc$, $\RT_q(\Th) \cap \Hdv$ (both for $q \geq 0$), and
similarly for the subspaces reflecting the different boundary conditions. The same notation will
also be used on the patches $\Ta$.

\subsection{\texorpdfstring{$\LT$}{L2}-orthogonal projectors and elementwise canonical \cred{projector}s}\label{sec_projs}

For $q \geq 0$, let $\Pi_h^q$ denote the $\Lti{\elm}$-orthogonal projector onto $\PP_q(\elm)$ or the elementwise $\Lt$-orthogonal projector onto \cred{the piecewise (broken) polynomials} $\PP_q(\Th)$, \ie, for $v \in \Lt$, $\Pi_h^q(v) \in \PP_q(\Th)$ is, separately for all $\elm \in \Th$, given by
\be \label{eq_L2proj}
    (\Pi_h^q (v), v_h)_\elm = (v, v_h)_\elm \qquad \forall v_h \in \PP_q(\elm).
\ee
Then, $\bPi_h^q$ is given componentwise by $\Pi_h^q$.

We will also use the $\Lt$-orthogonal projector $\RToproj{q}$ onto \cred{the broken Raviart--Thomas space} $\RT_q(\Th)$, given for $\tv \in \tLt$ also elementwise as: for all $\elm \in \Th$, $\RToproj{q}(\tv)\vert_\elm \in \RT_q(\elm)$ is such that
\bse\be \label{eq_RToproj}
    (\RToproj{q}(\tv),\tv_h)_\elm = (\tv,\tv_h)_\elm \qquad \forall \tv_h \in \RT_q(\elm),
\ee
or, equivalently,
\be \label{eq_RToproj_eq}
    \RToproj{q}(\tv)\vert_\elm \eq \arg \min_{\tv_h \in \RT_q(\elm)}
\norm{\tv - \tv_h}_\elm.
\ee\ese

\cred{Finally, we will also need the elementwise (broken) canonical projectors. Let $\tv \in \tLt$ such that $\tv\vert_\elm \in [C^1(\elm)]^3$ for all $\elm \in \Th$ be given. Below, we only use piecewise polynomial $\tv$ which satisfy these assumptions. Separately on each element $\elm \in \Th$, f}ollowing~\cite{Bof_Brez_For_MFEs_13, Ra_Tho_MFE_77}, the canonical $q$-degree Raviart--Thomas
\cred{projector} $\RTproj{q}(\tv)\cred{\vert_\elm} \in \RT_q(\elm)$, $q \geq 0$, is given by
\bse \label{eq_RT_proj} \bat{2}
\<\RTproj{q}(\tv) \scp \tn_\elm, r_h\>_\sd
&=
\<\tv \scp \tn_\elm, r_h\>_\sd
\qquad
& &
\forall r_h \in \PP_q(\sd), \quad \forall \sd \in \FK,
\label{eq_RT_proj_face}
\\
(\RTproj{q}(\tv), \tr_h)_\elm
&=
(\tv,\tr_h)_\elm
& &
\forall \tr_h \in [\PP_{q-1}(\elm)]^3,
\label{eq_RT_proj_vol}
\eat \ese
where $\tn_\elm$ is the unit outer normal vector of the element $\elm$.
\cred{Similarly, f}ollowing~\cite{Bof_Brez_For_MFEs_13, Ned_mix_R_3_80}, canonical $q$-degree N\'ed\'elec
\cred{projector} $\NDproj{q}(\tv)\cred{\vert_\elm} \in \ND_q(\elm)$, $q \geq 0$, is given,
\cred{separately on each element $\elm \in \Th$,} by
\bse \label{eq_ND_proj} \bat{2}
\<\NDproj{q}(\tv) \cdot \ttau_e,r_h \>_e
&
=
\< \tv \cdot \ttau_e, r_h \>_e
\qquad
& &
\forall r_h \in \PP_q(e), \quad \forall e \in \EK, \label{eq_ND_proj_1}
\\
\<\NDproj{q}(\tv) \times \tn_\elm,\tr_h \>_F
&=
\<\tv \times \tn_\elm,\tr_h \>_F
& &
\forall \tr_h \in [\PP_{q-1}(F)]^2, \quad \forall F \in \FK, \label{eq_ND_proj_2}
\\
(\NDproj{q}(\tv),\tr_h)_\elm
&=
(\tv,\tr_h)_\elm
& &
\forall \tr_h \in [\PP_{q-2}(\elm)]^3 \label{eq_ND_proj_3},
\eat \ese
where $\EK$ stands for the set of edges of $\elm$, $\ttau_e$ and $\<\cdot,\cdot\>_e$
respectively denote a (unit) tangential vector (the orientation does not matter)
and the $L^2(e)$ scalar product of the edge $e \in \EK$, and where we implicitly
\cred{complement $\tr_h$ by a zero component in the normal direction of the face $\sd$} in~\eqref{eq_ND_proj_2}.
These \cred{projector}s crucially satisfy, on \cred{each} tetrahedron $\elm$ and for all $\tv \in [C^1(\elm)]^3$, the commuting properties
\bse \label{eq_com_prop} \ba
    \Dv \RTproj{q}(\tv) & = \Pi_h^q (\Dv \tv), \label{eq_com_prop_div} \\
    \Crl (\NDproj{q} (\tv)) & = \RTproj{q}(\Crl \tv). \label{eq_com_prop_curl}
\ea \ese

\subsection{Functional inequalities} \label{sec_ineq}

We will need the four following functional inequalities.

First, from~\cite[Theorems~3.4 and~3.5]{Cost_Dau_Nic_sing_Maxw_99},
\cite[Theorem~2.1]{Hipt_Pechs_discr_reg_dec_19}, and the discussion
in~\cite[Section~3.2.1]{Chaum_Ern_Voh_Maxw_22}, it follows that for a
\cred{Lipschitz} polyhedral domain $\om \subset \Om$ with a \cred{Lipschitz}
Dirichlet boundary $\gD$ given by some of its faces, there exists
a constant $\Clift{\om}$ such that for all $\tv \in \HcDi{\om}$, there
exists $\ttw \in \tHoi{\om} \cap \HcDi{\om}$ such that $\Crl \ttw = \Crl \tv$ and
\be \label{eq_lift}
    \norm{\Gr \ttw}_\om \leq \Clift{\om} \norm{\Crl \tv}_\om.
\ee
When the $\gD$ has either a zero measure (in which case $\ttw$ can additionally be taken of mean
value zero componentwise) or coincides with $\pt \om$ and if $\om$ is convex, one can take
$\Clift{\om} = 1$, see~\cite{Cost_Dau_Nic_sing_Maxw_99} together with~\cite[Theorem~3.7]{Gir_Rav_NS_86}
for \cred{$\gD = \pt \om$} and~\cite[Theorem~3.9]{Gir_Rav_NS_86} for \cred{$\gD = \emptyset$}. Actually, $\ttw$ can even be taken in $\tHoDi{\om}$, though the above characterizations
$\Clift{\om} = 1$ may be lost. \cred{[--]}

Second, for any $\tv \in \tHoi{\om}$ of componentwise mean value zero on $\om$ or with the trace
equal to zero on $\gD \subset \pt \om$ which consists of at least one face of $\om$, the Poincar\'e--Friedrichs inequality gives
\be \label{eq_Poinc_Fried}
\norm{\tv}_{\om} \leq \CPFi{\om} h_\om \norm{\Gr \tv}_{\om}.
\ee
In the first case, $\CPFi{\om} \leq 1 / \pi$ for convex $\om$,
see~\cite{Pay_Wei_Poin_conv_60, Beben_Poin_conv_03}, whereas in the second case,
$\CPFi{\om} \leq 1$ when there exists a unit vector $\tm$ such that the straight
semi-line of direction $\tm$ originating at any point in $\om$ hits the boundary $\pt \om$
in the subset $\gD$, see~\cite{Voh_Poinc_disc_05}. We will employ~\eqref{eq_Poinc_Fried} with $\om=\oma$, where $\CPFi{\oma}$ only depends on the shape-regularity parameter of the mesh $\kappa_\Th$.

Third, for any mesh element $\elm \in \Th$, $\tv \in \tHoi{\elm}$, and $q \geq 0$, there holds the $hp$ approximation/Poincar\'e inequality
\be \label{eq_Poinc_hp}
\norm{\tv - \bPi_h^q(\tv)}_{\elm}
\leq
C_{\mathrm{hp}} \frac{h_\elm}{q+1} \norm{\Gr \tv}_{\elm}
\ee
for a generic constant $C_{\mathrm{hp}}$ only depending on $\kappa_\elm$
\cite{Bab_Sur_hp_FE_87}.

Finally,
the Poincar\'e--Friedrichs--Weber inequality, see~\cite[Proposition 7.4]{Fer_Gil_Maxw_BC_97}
and more precisely~\cite[Theorem A.1]{Chaum_Ern_Voh_Maxw_22} for the form of the constant,
will be useful: for all vertices $\ver \in \Vh$ and all vector-valued functions
$\tv \in \Hcdual{\oma} \cap \Hdva$ \cred{(recall the notation~\eqref{eq_loc_spaces_int}
or~\eqref{eq_loc_spaces_Dir})} with $\Dv \tv =0$, we have
\be
    \label{eq_weber_patch}
    \norm{\tv}_{\oma} \leq C_{\mathrm{PFW}} h_{\oma} \norm{\Crl \tv}_{\oma},
\ee
where $C_{\mathrm{PFW}}$ only depends on the mesh shape-regularity $\kappa_{\Th}$.
\cred{Note that $\oma$ is always homotopic to a ball and that if a Dirichlet
boundary is involved, it always corresponds to a simply-connected set.}
Strictly speaking, the inequality is established in~\cite[Theorem A.1]{Chaum_Ern_Voh_Maxw_22}
for edge patches, but the proof can be easily extended to vertex patches.

\section{Main results} \label{sec_main_res}

We present here our main results.

\subsection{A stable local commuting projector in \texorpdfstring{$\HDV$}{H(div)}} \label{sec_Phd}

We need to first recall \cred{(a slightly modified version of)} the projector $\Phd$ from~\cite[Definition~3.1]{Ern_Gud_Sme_Voh_loc_glob_div_22}:

\bd[A stable local commuting projector in $\HDV$] \label{def_Phd} Let $\ttw \in \HdvN$, \cred{a simplicial mesh $\Th$ of $\Om$, and a polynomial degree $p \geq 0$} be given.

\begin{enumerate} [parsep=0.2pt, itemsep=0.2pt, topsep=1pt, partopsep=1pt, leftmargin=13pt]

\item \label{dfd_1} Define a broken Raviart--Thomas polynomial $\btau_h \in \RT_p(\Th)$, on each mesh element~\cred{$\elm$}, via
\be \label{eq_tauh}
\btau_h\vert_\elm \eq \arg \min_{\substack{
\ttw_h \in \RT_p(\elm)
\\
\Dv \ttw_h = \Pi_h^p(\Dv \ttw)
}}
\norm{\ttw - \ttw_h}_\elm \qquad \forall \elm \in \Th.
\ee

\item \label{dfd_2} Define a Raviart--Thomas polynomial $\frh^\ver \in \RT_p(\Ta) \cap \Hdva$, on each vertex patch~\cred{$\Ta$}, via
\be \label{eq_sigma_a}
\frh^\ver
\eq
\arg \min_{\substack{
\ttw_h \in \RT_p(\Ta) \cap \Hdva
\\
\Dv \ttw_h= \Pi_h^p(\psia \Dv \ttw + \Gr \psia \scp \ttw)
}}
\norm[\big]{\RTproj{p} (\psia \btau_h) - \ttw_h}_\oma \qquad \forall \ver \in \Vh.
\ee

\item \label{dfd_3} Extending $\frh^\ver$ by zero outside of the patch subdomain $\oma$, define $\Phd (\ttw) \in \RT_p(\Th) \cap \HdvN$ via
\be \label{eq_sigma}
    \Phd (\ttw) \eq \frh \eq \suma \frh^\ver.
\ee

\end{enumerate}
\ed

The above definition works on the entire space $\HdvN$, in contrast to the canonical projector $\RTproj{p}$ from~\eqref{eq_RT_proj}. In Step~\ref{dfd_1}, we simply project $\ttw$ to the $\RT_p$ space elementwise, under the divergence constraint. The intermediate field $\btau_h$ is close to $\ttw$ but does not in general lie in $\HdvN$\cred{, possibly having a discontinuous normal trace}. This is corrected patchwise in Step~\ref{dfd_2}: a cut-off of $\btau_h$ is made by the hat function $\psia$, the canonical elementwise projector $\RTproj{p}$ from~\eqref{eq_RT_proj} is applied in order not to increase the polynomial degree by one, and a suitable divergence constraint, the elementwise $L^2$ projection~\eqref{eq_L2proj} of $\Dv(\psia \ttw)$ onto piecewise $p$-degree polynomials, is applied to crucially obtain the commuting property after Step~\ref{dfd_3}:
\be \label{eq_Phd_com} \bs
    \Dv \Phd (\ttw) & \cred{\reff{eq_sigma}= \Dv \frh = \suma \Dv \frh^\ver \reff{eq_sigma_a}= \suma \Pi_h^p(\Dv(\psia \ttw))} \\
    & \cred{= \Pi_h^p\bigg(\Dv\bigg(\suma (\psia \ttw)\bigg)\bigg)} = \Pi_h^p(\Dv \ttw),
\es \ee
\cred{using the linearity of the divergence operator, the linearity of the elementwise $L^2$ projection~\eqref{eq_L2proj}, and the partition of unity by the hat functions~\eqref{eq_PU}.}

\brem[Divergence constraint in~\eqref{eq_sigma_a}]\label{rem_contr_mod} The second term in the divergence constraint in~\eqref{eq_sigma_a} is modified with respect to~\cite[Definition~3.1]{Ern_Gud_Sme_Voh_loc_glob_div_22}, where $\Pi_h^p(\psia \Dv \ttw + \Gr \psia \scp \btau_h)$ is used instead. With this modification, a supplementary term arises on the right-hand side of equation~(4.8) in the proof of Lemma~4.6 in~\cite{Ern_Gud_Sme_Voh_loc_glob_div_22}. This term writes and can be bounded as
\be \label{eq_sup_term}
    (\Gr \psia \scp (\tv-\btau_{\mathcal{T}}), \Pi_{\mathcal{T}}^p \vf)_\oma \leq C \norm{\Gr \psia}_{\infty, \oma} h_\oma \norm{\tv-\btau_{\mathcal{T}}}_\oma \norm{\Gr \vf}_\oma \leq C \norm{\tv-\btau_{\mathcal{T}}}_\oma,
\ee
where the constant $C$ only depends on the shape regularity parameter $\kappa_{\Th}$, so that Lemma~4.6 and all results in~\cite{Ern_Gud_Sme_Voh_loc_glob_div_22} still hold. This subtle switch from $\btau_h$ to $\ttw$ in the second term in the divergence constraint in~\eqref{eq_sigma_a}, however, turns out crucial for our developments here, insuring in particular the well-posedness of Definition~\ref{def_Phc} below in that it makes~\eqref{eq_bha_min_set} in Lemma~\ref{lem_bha} hold true. \erem

\subsection{A stable local commuting projector in \texorpdfstring{$\HC$}{H(curl)}} \label{sec_Phc}

We define here our stable local commuting projector in $\HC$.

\bd[A locally defined mapping from $\HcN$ to $\ND_p(\Th) \cap \HcN$] \label{def_Phc} Let $\tv \in \HcN$, \cred{a simplicial mesh $\Th$ of $\Om$, and a polynomial degree $p \geq 0$} be given.

\begin{enumerate} [parsep=0.2pt, itemsep=0.2pt, topsep=1pt, partopsep=1pt, leftmargin=13pt]

\item \label{dfc_0} Set $\ttw \eq \Crl \tv$, so that $\ttw \in \HdvN$ with $\Dv \ttw =0$, and define $\btau_h \in \RT_p(\Th)$ by~\eqref{eq_tauh} and $\frh^\ver \in \RT_p(\Ta) \cap \Hdva$ by~\eqref{eq_sigma_a} from Definition~\ref{def_Phd}.

\item \label{dfc_1} Define a broken N\'ed\'elec polynomial $\biot_h \in \ND_p(\Th)$, on each mesh element~\cred{$\elm$}, via
\be \label{eq_bioth}
\biot_h\vert_\elm \eq \arg \min_{\substack{
\tv_h \in \ND_p(\elm)
\\
\Crl \tv_h = \btau_h
}}
\norm{\tv - \tv_h}_\elm \qquad \forall \elm \in \Th.
\ee

\item \label{dfc_2} Define a Raviart--Thomas polynomial $\btha \in \RT_{p+1}(\Ta) \cap \Hdva$, on each vertex patch~\cred{$\Ta$}, via
\be \label{eq_tha}
\btha \eq \arg \hspace*{-1cm}
\min_{\substack{
\tv_h \in \RT_{p+1}(\Ta) \cap \Hdva
\\
\Dv \tv_h= \Pi_h^{p+1}(- \Gr \psia \scp (\Crl \tv))
\\
(\tv_h, \tr_h)_\elm = (\Gr \psia \vp \biot_h, \tr_h)_\elm
\quad
\forall \tr_h \in [\PP_0(\elm)]^{3}, \, \forall \elm \in \Ta}} \hspace*{-1cm} \norm{\Gr \psia \vp \biot_h- \tv_h}_\oma \qquad \forall \ver \in \Vh.
\ee

\item \label{dfc_3} Extending $\btha$ by zero outside of the patch subdomain $\oma$, define
\be \label{eq_td}
    \bdl_h \eq \suma \btha,
\ee
which gives
\be \label{eq_td_prop}
    \bdl_h \in \RT_{p+1}(\Th) \cap \HdvN \quad \text{ and } \quad \Dv \bdl_h = 0.
\ee

\item \label{dfc_4} For all tetrahedra $\elm \in \Th$, consider the $(p+1)$-degree Raviart--Thomas elementwise minimizations
\ba \label{eq_tda}
    \bdla\vert_\elm \eq \arg \hspace*{-0.65cm} \min_{\substack{\tv_h \in \RT_{p+1}(\elm)\\ \Dv \tv_h= 0\\ \tv_h \scp \tn_\elm = \RTproj{p+1}(\psia \bdl_h) \scp \tn_\elm \, \text{ on } \pt \elm }} \hspace*{-0.65cm} \norm[\big]{\RTproj{p+1}(\psia \bdl_h) - \tv_h}_\elm \qquad \forall \ver \in \VK,
\ea
which gives
\bse \label{eq_div_free_dec} \ba
    \bdla & \in \RT_{p+1}(\Ta) \cap \Hdva \quad \text{ and } \quad \Dv \bdla= 0 \qquad \forall \ver \in \Vh, \label{eq_bdla_div} \\
    \bdl_h & = \suma \bdla. \label{eq_bdla_sum}
\ea \ese

\item \label{dfc_5} Define a N\'ed\'elec polynomial $\bha \in \ND_p(\Ta) \cap \Hca$, on each vertex patch~\cred{$\Ta$}, via
\be \label{eq_bha}
\bha
\eq
\arg \min_{\substack{
\tv_h \in \ND_p(\Ta) \cap \Hca
\\
\Crl \tv_h = \frh^\ver + \RTproj{p} (\btha - \bdla)
}}
\norm[\big]{\NDproj{p} (\psia \biot_h) - \tv_h}_\oma \qquad \forall \ver \in \Vh.
\ee

\item \label{dfc_6} Extending $\bha$ by zero outside of the patch subdomain $\oma$, define
\be \label{eq_bh}
    \Phc (\tv) \eq \bh_h \eq \suma \bha,
\ee
which gives
\be \label{eq_bh_prop}
    \Phc (\tv) \in \ND_p(\Th) \cap \HcN.
\ee

\end{enumerate}

\ed

It is not obvious from Definition~\ref{def_Phc} that the operator $\Phc$ is well defined. This is justified in Section~\ref{sec_Phc_WP} below. We, on the other hand, immediately see:

\brem[Definition~\ref{def_Phc}] \label{rem_Phc} The operator~$\Phc$ is defined over the entire infinite-dimensional space $\HcN$, is composed of simple piecewise polynomial projections, and is defined locally, in a neighborhood of each mesh vertex. \cred{The operator~$\Phc$ is, however, not defined over the entire $\tLt$, since on Step~\ref{dfc_0} of Definition~\ref{def_Phc}, we need $\ttw \eq \Crl \tv$ to be well defined, to lie in $\HdvN$, and to satisfy $\Dv \ttw =0$.} \erem

The purpose of Step~\ref{dfc_0} is to prepare suitable piecewise polynomial data approximating $\Crl \tv$ following Definition~\ref{def_Phd}, which will also later ensure the commuting property with the operator $\Phd$.

The purpose of Step~\ref{dfc_1} is to bring all further considerations from the infinite-dimensional level of the given function $\tv \in \HcN$ to the finite-dimensional level of $\biot_h \in \ND_p(\Th)$. In particular, all subsequent steps only involve $\biot_h$ and other piecewise polynomials, which in particular justifies the use of the canonical \cred{elementwise} Raviart--Thomas and N\'ed\'elec \cred{projector}s $\RTproj{}$ and $\NDproj{}$ \cred{from~\eqref{eq_RT_proj} and~\eqref{eq_ND_proj}} in Steps~\ref{dfc_4} and~\ref{dfc_5}.

At the continuous level, $\Gr \psia \vp \tv$ would belong to the $\Hdva$ space, with the divergence being equal to $- \Gr \psia \scp (\Crl \tv)$, \cf\ \cite[equations~(4.4)--(4.\cred{6})]{Chaum_Voh_Maxwell_equil_23}. In Step~\ref{dfc_2}, we mimic this at the discrete level with $\btha$. Contrarily to the continuous level, where $\suma (\Gr \psia \vp \tv) = \nv$, on Step~\ref{dfc_3}, $\bdl_h = \suma \btha \neq \nv$, though one may hope that $\bdl_h \approx \nv$. One, in turn, immediately sees that $\Dv \bdl_h = \suma \Dv \btha = 0$ by the partition of unity~\eqref{eq_PU}, so that $\bdl_h$ is a divergence-free Raviart--Thomas polynomial.

The purpose of Step~\ref{dfc_4} is to achieve the decomposition of $\bdl_h$ by the divergence-free local contributions $\bdla$ as per~\eqref{eq_div_free_dec}. This procedure is taken from~\cite[Section~5.1]{Chaum_Voh_Maxwell_equil_23} and is crucial here. Its key ingredient is the (non-traditional) additional constraint in~\eqref{eq_tha} on orthogonality with respect to piecewise vector-valued constants.

Finally, as in the $\HC$-equilibration of~\cite[Section~5.2]{Chaum_Voh_Maxwell_equil_23}, all this preparatory work allows us to pose the local problem~\eqref{eq_bha} in Step~\ref{dfc_5}; in particular, the requested curl given by $\frh^\ver + \RTproj{p} (\btha - \bdla)$ has to belong to $\RT_p(\Ta) \cap \Hdva$ and be divergence-free. In other words, Steps~\ref{dfc_2}--\ref{dfc_4} serve to prepare this crucial datum for problem~\eqref{eq_bha}; as for Steps~\ref{dfc_1}, \ref{dfc_5}, and~\ref{dfc_6}, these fully mimic Steps~\ref{dfd_1}--\ref{dfd_3} from Definition~\ref{def_Phd}. In particular, the definition of $\Phc (\tv)$ can finally be finished in Step~\ref{dfc_6}.

\brem[Constraints in~\eqref{eq_tauh} and~\eqref{eq_bioth}] The divergence and curl constraints in respectively~\eqref{eq_tauh} and~\eqref{eq_bioth} are not essential for the projectors $\Phd$ and $\Phc$ and could be removed. Then, actually, the modification of the divergence constraint in~\eqref{eq_sigma_a} of Remark~\ref{rem_contr_mod} would not be necessary, since then the use of $\btau_h$ or $\ttw$ in~\eqref{eq_sigma_a} would be equivalent.
However, the use of constraints in~\eqref{eq_tauh} and~\eqref{eq_bioth} leads to a slightly sharper constants in~\eqref{eq_approx}--\eqref{eq_stab_curl_norm} and is \cred{important} for the $hp$-approximation result of Theorem~\ref{thm_hp}. \erem

Recall the $\Lt$-orthogonal projector $\RToproj{q}$ from~\eqref{eq_RToproj}. Our first main result is:

\bt[Commutativity, projection, approximation, and stability of $\Phc$] \label{thm_Phc} Let a \cred{simplicial} mesh $\Th$ of $\Om$ and a polynomial degree $p \geq 0$ be fixed. Then, the operator~$\Phc$ from Definition~\ref{def_Phc} is a commuting projector since
\bat{2}
\Crl \Phc (\tv) & = \Phd (\Crl \tv) \qquad \qquad  & & \forall \tv \in \HcN, \label{eq_com} \\
\Phc (\tv) & = \tv & & \forall \tv \in \ND_p(\Th) \cap \HcN. \label{eq_proj}
\eat
Moreover, $\Phc$ has optimal approximation properties of an elementwise unconstrained $\tLT$-orthogonal projector and is stable in that for any function $\tv \in \HcN$ and any mesh element $\elm \in \Th$, there holds
\ba
\norm[\big]{\tv - \Phc (\tv)}_\elm^2 \quad {} & \nn \\
+ \bigg(\frac{h_\elm}{p+1}\norm[\big]{\Crl\big(\tv-\Phc (\tv)\big)}_\elm\bigg)^2 \ls {} &  \sum_{\elmtt\in\TK} \Bigg\{\min_{\tv_h \in \ND_p(\elmtt)}\norm{\tv - \tv_h}_\elmtt^2 \nn\\
{} & + \bigg(\frac{h_{\elmtt}}{p+1}\norm{\Crl\tv - \RToproj{p}(\Crl\tv)}_\elmtt\bigg)^2\Bigg\},\label{eq_approx} \\
\norm[\big]{\Phc (\tv)}_\elm^2 \ls {} & \sum_{\elmtt \in \TK} \Bigg\{\norm{\tv}_\elmtt^2 \nn \\ {} & + \Bigg(\frac{h_\elmtt}{p+1} \norm[\big]{\Crl \tv - \RToproj{p} (\Crl \tv)}_\elmtt\Bigg)^2\Bigg\},\label{eq_stab} \\
\norm[\big]{\Phc (\tv)}_\elm^2 + h_\Om^2\norm[\big]{\Crl \Phc (\tv)}_\elm^2 \ls {} & \sum_{\elmtt \in \TK}\big\{\norm{\tv}_{\elmtt}^2 + h_\Om^2 \norm{\Crl \tv}_\elmtt^2\big\}, \label{eq_stab_curl_norm}
\ea
where $\TK$ collects the elements $\elmtt$ of $\Th$ sharing a vertex with $\elm$ or its neighbor and $h_\Om$ denotes the diameter of $\Omega$. The constant hidden in $\ls$ only depends on the shape-regularity parameter $\kappa_{\Th}$ of the mesh $\Th$ and the polynomial degree $p$.%
\et

We will prove Theorem~\ref{thm_Phc} in Section~\ref{sec_Phc_proof} below. As stated, it shows that $\Phc$ satisfies all the properties~\ref{pr_ent_sp}--\ref{pr_proj} of Section~\ref{sec_intr}:

\brem[Theorem~\ref{thm_Phc}] We remark that~\eqref{eq_com} is precisely the commuting property desired in the middle column of~\eqref{eq_ex_seq}, whereas~\eqref{eq_proj} is the projector property. Moreover, \eqref{eq_approx} is the optimal approximation property:
the first term on the right-hand side of~\eqref{eq_approx} is the local-best (elementwise) approximation error, that of the $\tLT$-orthogonal projector onto $\ND_p(\elmtt)$, without any constraint on the curl, whereas the second term on the right-hand side of~\eqref{eq_approx},
\be \label{eq_hp_data_osc}
    \cred{\frac{h_\elmtt}{p+1} \norm[\big]{\Crl \tv - \RToproj{p} (\Crl \tv)}_\elmtt,}
\ee
is an $hp$ ``data oscillation term'' \cred{known from a posteriori error analysis}, which in particular disappears when $\Crl \tv$ is a \cred{piecewise $p$-degree Raviart--Thomas polynomial, $\Crl \tv \in \RT_p(\Th) \cap \HdvN$, or, equivalently, since $\Dv (\Crl \tv) = 0$, $\Crl \tv \in [\PP_p(\Th)]^3 \cap \HdvN$}. Finally, \eqref{eq_stab} is stability in the $\tLt$-norm, up to the $hp$ data oscillation \cred{terms~\eqref{eq_hp_data_osc}, \ie, stability in the $\tLt$-norm when the terms~\eqref{eq_hp_data_osc} are zero and stability in the $\tLt$-norm up to the terms~\eqref{eq_hp_data_osc}, which will typically be (much) smaller than $\norm{\tv}_\elmtt$ itself. This is in coincidence with the fact that our projector is not defined over the entire $\tLt$, see Remark~\ref{rem_Phc}. Finally,} \eqref{eq_stab_curl_norm} is stability in the $\Hc$-norm, where the (physical) scaling by $h_\Om$ is chosen for the curl term (other scalings, by at least the local mesh sizes $h_\elm$, could be chosen as well). \erem

\cred{\brem [Statements~\eqref{eq_approx} and~\eqref{eq_stab} with $h_\elmtt$ versus $h_\elmtt/(p+1)$] In Theorem~\ref{thm_Phc}, it may appear strange to employ the $hp$ terms $h_\elmtt/(p+1)$ in place of $h_\elmtt$ (and similarly with $h_\elm$) since anyhow the hidden constant depends on the polynomial degree $p$. The statements hold equivalently if $h_\elmtt/(p+1)$ are replaced by $h_\elmtt$ (and similarly with $h_\elm$). We use this formulation because 1) it matches with our presentation of Theorem~\ref{thm_hp} below (where constants are independent of $p$ and such a writing is crucial); and 2) it is convenient for us in the proof of Theorem~\ref{thm_hp} on pages~\pageref{page_proof_Thm_3_beg}--\pageref{page_proof_Thm_3_end}. There it becomes crucial to use the current (re)writing. \erem}

\subsection{Equivalence of local-best and global-best approximations in \texorpdfstring{$\HC$}{H(curl)} \cred{on general domains and} including data oscillation} \label{sec_equiv_LG}

The results of Theorem~\ref{thm_Phc} immediately lead to the following extensions of~\cite[Theorem~2]{Chaum_Voh_loc_glob_curl_21}:

\bt[Equivalence of local-best and global-best approximations in $\HC$ \cred{on general domains and} including data oscillation] \label{thm_loc_glob} Let $\tv \in \HcN$, a \cred{simplicial} mesh $\Th$ of $\Om$, and a polynomial degree $p \geq 0$ be fixed. Then
\bse \be \label{eq_LG_constr} \bs
{} & \min_{\substack{\tv_h \in \ND_p(\Th) \cap \HcN\\\Crl\tv_h = \Phd (\Crl \tv)}}\norm{\tv - \tv_h}^2
+ \sum_{\elm \in \Th} \bigg(\frac{h_{\elm}}{p+1}\norm{\Crl\tv - \RToproj{p}(\Crl\tv)}_\elm\bigg)^2 \\
\approx {} & \sum_{\elm\in\Th} \bigg[\min_{\tv_h \in \ND_p(\elm)}\norm{\tv - \tv_h}_\elm^2 + \bigg(\frac{h_{\elm}}{p+1}\norm{\Crl\tv - \RToproj{p}(\Crl\tv)}_\elm\bigg)^2\bigg],
\es \ee
and
\be \label{eq_LG_unconstr} \bs
{} & \min_{\tv_h \in \ND_p(\Th) \cap \HcN}\bigg[\norm{\tv - \tv_h}^2 + \sum_{\elm \in \Th}\bigg(\frac{h_{\elm}}{p+1}\norm{\Crl(\tv-\tv_h)}_\elm\bigg)^2\bigg] \\
\approx {} & \sum_{\elm\in\Th} \bigg[\min_{\tv_h \in \ND_p(\elm)}\norm{\tv - \tv_h}_\elm^2 + \bigg(\frac{h_{\elm}}{p+1}\norm{\Crl\tv - \RToproj{p}(\Crl\tv)}_\elm\bigg)^2\bigg],
\es \ee \ese
where the hidden constants only depend on the shape-regularity parameter $\kappa_{\Th}$ of the mesh $\Th$ and the polynomial degree $p$. \et

The constraint in~\eqref{eq_LG_constr} uses the projector $\Phd$ from Definition~\ref{def_Phd}. By the projection property of $\Phd$ on $\RT_p(\Th) \cap \HdvN$, there immediately follows that $\Phd(\Crl\tv) = \Crl\tv$ when $\Crl \tv \in [\PP_p(\Th)]^3$, since in this case $\Crl \tv \in \RT_p(\Th) \cap \HdvN$.
Thus~\eqref{eq_LG_constr} indeed extends~\cite[Theorem~2]{Chaum_Voh_loc_glob_curl_21} to the case when the datum $\tv \in \HcN$ has a non-polynomial curl, $\Crl \tv \not \in [\PP_p(\Th)]^3$. \cred{Moreover, the proofs in~\cite{Chaum_Voh_loc_glob_curl_21} require the sequence~\eqref{eq_ex_seq} to be exact (which happens when $\Omega$ is homotopic to a ball (contractible) and either $\GD = \pt \Om$ or $\GN = \pt \Om$), which we do not suppose here.} Remarkably, both the tangential trace continuity and the curl constraint from the left-hand side of~\eqref{eq_LG_constr} are removed in the right-hand side of~\eqref{eq_LG_constr}.
Finally, in~\eqref{eq_LG_unconstr}, only the tangential trace continuity constraint is removed,
but the simultaneous approximation of the curl is the best-possible.

\begin{remark}[Mixed finite element discretization]
The result~\eqref{eq_LG_constr} includes the projector $\Phd$ in its constraint,
which may appear as an arbitrary choice. In contrast, the following estimates may have a more direct application: if $p \geq 1$, $\Omega$ is convex, either $\GD = \emptyset$ or $\GN = \emptyset$, and the mesh $\Th$ is quasi-uniform in that $h_\elm \approx h \eq \max_{\elmtt \in \Th} h_\elmtt$ for all $\elm \in \Th$, we have
\be \label{eq_LG_constr_glob} \bs
{} & \min_{\substack{\tv_h \in \ND_p(\Th) \cap \HcN\\\Crl\tv_h = \bPi_h^{p,\dv}(\Crl \tv)}}\norm{\tv - \tv_h}^2
+ \sum_{\elm \in \Th} \bigg(\frac{h_{\elm}}{p+1}\norm{\Crl\tv - \RToproj{p}(\Crl\tv)}_\elm\bigg)^2 \\
\approx {} & \sum_{\elm\in\Th} \bigg[\min_{\tv_h \in \ND_p(\elm)}\norm{\tv - \tv_h}_\elm^2 + \bigg(\frac{h_{\elm}}{p+1}\norm{\Crl\tv - \RToproj{p}(\Crl\tv)}_\elm\bigg)^2\bigg],
\es \ee
where the hidden constant only depends on the shape-regularity parameter $\kappa_{\Th}$ of the mesh $\Th$ and the polynomial degree $p$ and where
\be \label{eq_RT_glob}
\bPi_h^{p,\dv}(\Crl \tv)\eq
\arg \min_{\substack{
\ttw_h \in \RT_p(\Th) \cap \HdvN
\\
\Dv \ttw_h= 0
}}
\norm[\big]{\Crl \tv - \ttw_h}^2
\ee
is the global $\tLt$-orthogonal projection of $\Crl \tv$ onto the divergence-free subspace of the Raviart--Thomas space
$\RT_p(\Th) \cap \HdvN$. Indeed, the constrained minimization problem~\eqref{eq_RT_glob}
naturally arises in the mixed finite element discretization of the curl--curl problem,
consisting in finding
\[
    \tv_h \eq \arg \min_{\substack{\tv_h \in \ND_p(\Th) \cap \HcN\\\Crl\tv_h = \bPi_h^{p,\dv} (\Crl \tv)}}\norm{\tv - \tv_h}^2,
\]
\ie, $\tv_h \in \ND_p(\Th) \cap \bm{H}_{0,\mathrm{N}}(\crl, \Om)$,
$\tp_h \in \RT_p(\Th) \cap \HdvN$, and $q_h \in \PP_p(\Th)$, of mean value zero when $\GD = \emptyset$, such that
\begin{equation*}
\begin{array}{cccccccll}
(\tv_h,\ttw_h) &\!+\!& (\tp_h,\Crl \ttw_h) & & &\!=\!& (\tv,\ttw_h) & \forall \ttw_h \in \ND_p(\Th) \cap \bm{H}_{0,\mathrm{N}}(\crl, \Om),
\\
(\Crl \tv_h,\tr_h) & &               &\!+\!& (q_h,\Dv \tr_h) &\!=\!& (\Crl \tv,\tr_h) & \forall \tr_h \in \RT_p(\Th) \cap \HdvN,
\\
              & &  (\Dv \tp_h,s_h)   & &                 &\!=\!& 0 & \forall s_h \in \PP_p(\Th).
\end{array}
\end{equation*}
Unfortunately, our current proof of~\eqref{eq_LG_constr_glob} only holds under the above fairly restrictive assumptions
on the domain and the mesh, and it remains an open question whether
\eqref{eq_LG_constr_glob} holds in more general settings.
We also remark that~\eqref{eq_LG_constr} and~\eqref{eq_LG_constr_glob} are
equivalents of~\cite[Theorem~3.3]{Ern_Gud_Sme_Voh_loc_glob_div_22} in the $\HcN$ context.
\end{remark}

We will prove Theorem~\ref{thm_loc_glob} and property~\eqref{eq_LG_constr_glob} in
Section~\ref{sec_loc_glob_proof} below.

\subsection{Optimal $hp$ approximation estimates in \texorpdfstring{$\HC$}{H(curl)}} \label{sec_hp}

We finally introduce an $hp$ approximation estimate, an equivalent
of~\cite[Theorem~3.6]{Ern_Gud_Sme_Voh_loc_glob_div_22} in the $\HC$-setting. It is optimal with respect to both the mesh size $h$ and the polynomial degree $p$
for arbitrary Sobolev regularity indices $s \geq 0$ for the approximated function
and $t \geq 0$ for its curl, where, importantly, these additional regularities are only
requested elementwise.

\bt[$hp$-optimal approximation estimates in $\HC$ under minimal Sobolev regularity]
\label{thm_hp}
Let $\tv \in \HcN$, a \cred{simplicial} mesh $\Th$ of $\Om$, and a polynomial degree $p \geq 0$ be fixed. Let
\be \label{eq_reg_v}
    \tv\vert_\elm \in {\bm H}^s(\elm) \qquad \forall \elm \in \Th
\ee
for a fixed regularity exponent $s \geq 0$.

\cred{We consider three cases.

\noindent {\em Case~(i) (piecewise polynomial curl).} Let $\Crl \tv \in [\PP_{p-1}(\Th)]^3$. Then
\be \label{eq_hp_pol_curl} \bs
{} & \min_{\tv_h \in \ND_p(\Th) \cap \HcN}\bigg[\norm{\tv - \tv_h}^2 + \sum_{\elm \in \Th}\bigg(\frac{h_{\elm}}{p+1}\norm{\Crl(\tv-\tv_h)}_\elm\bigg)^2\bigg] \\
\ls {} & \sum_{\elm\in\Th} \bigg[\bigg(\frac{h_\elm^{\min\{p+1,s\}}}{(p+1)^s} \norm{\tv}_{{\bm H}^s(\elm)}\bigg)^2\bigg],
\es \ee
where the hidden constant only depends on the shape-regularity parameter $\kappa_{\Th}$ of the mesh $\Th$ and the regularity exponent $s$.}

\noindent \cred{{\em Case~(ii) (convex patches and tangential boundary conditions)}.} Let the patch subdomains $\oma$ be convex for all vertices $\ver \in \Vh$ and let $\GD = \emptyset$. Let
\be \label{eq_reg_curl_v}
    (\Crl \tv)\vert_\elm \in {\bm H}^t(\elm) \qquad \forall \elm \in \Th
\ee
\cred{for a fixed regularity exponent $t \geq 0$}. Then
\be \label{eq_hp} \bs
{} & \min_{\tv_h \in \ND_p(\Th) \cap \HcN}\bigg[\norm{\tv - \tv_h}^2 + \sum_{\elm \in \Th}\bigg(\frac{h_{\elm}}{p+1}\norm{\Crl(\tv-\tv_h)}_\elm\bigg)^2\bigg] \\
\ls {} & \sum_{\elm\in\Th} \bigg[\bigg(\frac{h_\elm^{\min\{p+1,s\}}}{(p+1)^s} \norm{\tv}_{{\bm H}^s(\elm)}\bigg)^2 + \bigg(\frac{h_\elm}{p+1} \frac{h_\elm^{\min\{p+1,t\}}}{(p+1)^t} \norm{\Crl \tv}_{{\bm H}^t(\elm)} \bigg)^2\bigg],
\es \ee
where the hidden constant only depends on the shape-regularity parameter $\kappa_{\Th}$ of the mesh $\Th$ and the regularity exponents $s$ and $t$.

\noindent \cred{{\em Case~(iii) (general case)} Let~\eqref{eq_reg_curl_v} hold for a fixed regularity exponent $t \geq 0$. Then
\be \label{eq_hp_gen} \bs
{} & \min_{\tv_h \in \ND_p(\Th) \cap \HcN}\bigg[\norm{\tv - \tv_h}^2 + \sum_{\elm \in \Th}\bigg(h_{\elm}\norm{\Crl(\tv-\tv_h)}_\elm\bigg)^2\bigg] \\
\ls {} & \sum_{\elm\in\Th} \bigg[\bigg(\frac{h_\elm^{\min\{p+1,s\}}}{(p+1)^s} \norm{\tv}_{{\bm H}^s(\elm)}\bigg)^2 + \bigg(h_{\elm} \frac{h_\elm^{\min\{p+1,t\}}}{(p+1)^t} \norm{\Crl \tv}_{{\bm H}^t(\elm)} \bigg)^2\bigg],
\es \ee
where the hidden constant only depends on the shape-regularity parameter $\kappa_{\Th}$ of the mesh $\Th$ and the regularity exponents $s$ and $t$.}
\et

\cred{
\brem[The three cases in Theorem~\ref{thm_hp}]
Case~(i) in Theorem~\ref{thm_hp} describes the ``no data oscillation'' situation.
The requirement on convex patch subdomains $\oma$ in case~(ii) may be satisfied on
convex domains even for solutions with low regularity (e.g. stemming from discontinuous
coefficients), but it excludes nonconvex domains. From case~(iii), however, we see that
we lose at worse the factor $1/(p+1)$, and this only in the curl term.
\erem
}

\section{Well-posedness of Definition~\ref{def_Phc} of \texorpdfstring{$\Phc$}{$P_h^{p,\crl}$} } \label{sec_Phc_WP}

In this section, we establish that the operator $\Phc$ from Definition~\ref{def_Phc} is well defined. Let $\tv \in \HcN$ be fixed.

\bl[Problem~\eqref{eq_bioth}] \label{lem_bioth} For each mesh element $\elm \in \Th$, problem~\eqref{eq_bioth} for $\biot_h\vert_\elm$ is well defined. It can equivalently be stated by its Euler--Lagrange optimality conditions which read: find $\biot_h \in \ND_p(\elm)$ with $\Crl \biot_h = \btau_h$ such that
\be \label{eq_bioth_EL}
    (\biot_h, \tv_h)_\elm = (\tv, \tv_h)_\elm \qquad \forall \tv_h \in \ND_p(\elm) \text{ with } \Crl \tv_h = \nv.
\ee
\el

\bp As we take $\ttw \eq \Crl \tv$ in~\eqref{eq_tauh}, we have $\Dv \ttw = \Dv (\Crl \tv) = 0$, so that $\btau_h \in \RT_p(\elm)$ from~\eqref{eq_tauh} is divergence-free. Consequently, the minimization set in~\eqref{eq_bioth} is nonempty, \cf, \eg, \cite[equation~(2.3.62)]{Bof_Brez_For_MFEs_13} which gives
\be \label{eq_curl_Np_map_loc}
    \Crl \ND_p(\elm) = \{\tv_h \in \RT_p(\elm); \, \Dv \tv_h = 0\}.
\ee
Consequently, the convex minimization~\eqref{eq_bioth} is well posed and its Euler--Lagrange optimality conditions read as~\eqref{eq_bioth_EL}.
\ep

\bl[Problem~\eqref{eq_tha}] \label{lem_tha} For each mesh vertex $\ver \in \Vh$, problem~\eqref{eq_tha} for $\btha$ is well defined. Moreover, \eqref{eq_td_prop} holds. \el

\bp Problem~\eqref{eq_tha} is, on a first sight, over-constrained; the orthogonality with respect to piecewise vector-valued constants $\tr_h$ adds additional constraints to a well-posed problem without it.
It, however, by construction fits the framework of~\cite[Theorem~A.2]{Chaum_Voh_Maxwell_equil_23}, which ensures existence and uniqueness of $\btha$. Indeed, let $q \eq q' \eq p+1$, $g^\ver \eq (- \Gr \psia \scp (\Crl \tv))$, and $\btau_h^\ver \eq \Gr \psia \vp \biot_h$. Then, we only need to verify Assumption~A.1 of~\cite{Chaum_Voh_Maxwell_equil_23}, which requests
\bse\ba
g^\ver & \in \Lti{\oma} \quad \text{ and } \quad \btau_h^\ver \in \RT_{p+1}(\Ta), \label{eq_g_tau}
\\
(g^\ver, 1)_\oma & = 0 \quad \text{ when } \ver \not \in \overline \GD, \label{eq_ga_comp}
\\
(\btau_h^\ver, \Gr q_h)_\oma + (g^\ver, q_h)_\oma & = 0
\qquad
\forall q_h \in \PP_1(\Ta)\cap\Hsa. \label{eq_orth_A}
\ea \ese

Requirement~\eqref{eq_g_tau} is obvious. As for~\eqref{eq_ga_comp}, we have
\[
    (- \Gr \psia \scp (\Crl \tv),1)_\oma = -(\Crl \tv,\Gr \psia)_\oma = (\underbrace{\Dv (\Crl \tv)}_{=0}, \psia)_\oma - \underbrace{\<(\Crl \tv) \scp \tn, \psia\>_{\pt \oma}}_{=0}
\]
by the Green theorem for a vertex $\ver$ not being part of $\overline \GD$. Finally, \cred{for~\eqref{eq_orth_A},} let $q_h \in \PP_1(\Ta)\cap\Hsa$. Then, again by the Green theorem,
\[
    - (\Gr \psia \scp (\Crl \tv), q_h)_\oma = (\Dv(\Gr \psia \vp \tv), q_h)_\oma = - (\Gr \psia \vp \tv, \Gr q_h)_\oma.
\]
Consequently, since $(\cred{\td} \vp \cred{\te}) \scp \cred{\tf} = \cred{\te} \scp (\cred{\tf} \vp \cred{\td})$ for vectors $\cred{\td}, \cred{\te}, \cred{\tf} \in \RR^3$,
\ban
    {} & (\Gr \psia \vp \biot_h, \Gr q_h)_\oma - (\Gr \psia \scp (\Crl \tv), q_h)_\oma \\
    = {} & \cred{\sum_{\elm \in \Ta}}(\biot_h - \tv, \hspace*{-0.3cm}\underbrace{\Gr q_h \vp \Gr \psia}_{\in \ND_p(\elm) \text{ with } \Crl (\Gr q_h \vp \Gr \psia) = \nv}\hspace*{-0.3cm})_\elm \reff{eq_bioth_EL}= 0
\ean
from the Euler--Lagrange conditions~\eqref{eq_bioth_EL}.

As for~\eqref{eq_td_prop}, $\bdl_h \in \RT_{p+1}(\Th) \cap \HdvN$ is obvious since the contributions $\btha$, after extensions by zero outside of the patch subdomains $\oma$, belong to $\RT_{p+1}(\Th)$ and are $\HdvN$-conforming. The divergence-free property is then immediately seen from the partition of unity~\eqref{eq_PU} which gives
\[
    \Dv \bdl_h = \suma \Dv \btha = \suma [\Pi_h^{p+1}(- \Gr \psia \scp (\Crl \tv))] = 0.
\]
\ep

\bl[Problem~\eqref{eq_tda}] \label{lem_delta} For all tetrahedra $\elm \in \Th$ and all vertices $\ver \in \VK$, problem~\eqref{eq_tda} for $\bdla\vert_\elm$ is well defined.
Moreover, properties~\eqref{eq_div_free_dec} hold. \el

\bp

Let $\elm \in \Th$ be fixed. \cred{Problem~\eqref{eq_tda} is a Neumann problem which is well-posed when the Neumann compatibility condition
\[
    \< \RTproj{p+1}(\psia \bdl_h) \scp \tn_\elm, 1 \>_{\pt \elm} = 0
\]
holds. We affirm that this is the case when}
\be \label{eq_bdl_orth}
    (\bdl_h, \tr_h)_\elm = 0 \qquad \forall \tr_h \in [\PP_0(\elm)]^3
\ee
is satisfied. \cred{Indeed, using the property~\eqref{eq_RT_proj_face} of the canonical Raviart--Thomas projector, the Green theorem, the fact that $\Dv \bdl_h = 0$, and~\eqref{eq_bdl_orth}, we see
\ban
    \< \RTproj{p+1}(\psia \bdl_h) \scp \tn_\elm, 1 \>_{\pt \elm} & \refvs{\eqref{eq_RT_proj_face}}{=}{xxxx} \< (\psia \bdl_h) \scp \tn_\elm, 1 \>_{\pt \elm} = \< \bdl_h \scp \tn_\elm, \psia \>_{\pt \elm} \\
    & \refvs{}{=}{xxxx} (\Dv \bdl_h, \psia)_\elm + (\bdl_h, \Gr \psia)_\elm \refff{eq_td_prop}{eq_bdl_orth}= 0.
\ean
P}roperty~\eqref{eq_bdl_orth} is, \cred{in turn,} a simple consequence of the orthogonality assumption in~\eqref{eq_tha}, of definition~\eqref{eq_td}, and of the partition of unity~\eqref{eq_PU}, which altogether give
\[
    (\bdl_h, \tr_h)_\elm \reff{eq_td}= \sum_{\ver \in \VK} (\btha, \tr_h)_\elm \reff{eq_PU}= \sum_{\ver \in \VK} (\btha - \Gr \psia \vp \biot_h, \tr_h)_\elm \reff{eq_tha}= 0.
\]

\cred{Property~\eqref{eq_bdla_div} is an immediate consequence of the definition~\eqref{eq_tda}. The decomposition property}~\eqref{eq_bdla_sum} is then established in~\cite[Theorem~B.1, \cred{equation~(B.5)}]{Chaum_Voh_Maxwell_equil_23} (for the choice $q=q'=p+1$). \cred{It relies on the fact that the normal traces on faces of $\bdl_h$ and $\suma \bdla$ are equal from the imposition of $\bdla \scp \tn_\elm$ in~\eqref{eq_tda}, using the partition of unity~\eqref{eq_PU} together with the linearity and projection property of $\RTproj{}$ in~\eqref{eq_RT_proj}, whereas the interior moments vanish from the orthogonality requirement behind~\eqref{eq_tda}.}
\ep

\bl[Problem~\eqref{eq_bha}] \label{lem_bha} For each mesh vertex $\ver \in \Vh$, problem~\eqref{eq_bha} is well defined. Moreover, \eqref{eq_bh_prop} holds.  \el

\bp
The convex minimization~\eqref{eq_bha} is well posed if the minimization set in~\eqref{eq_bha}
is nonempty.

\cred{Recall that we suppose $\oma$ homotopic to a ball (contractible),
and that the essential boundary conditions in $\Hca$ and $\Hdva$ are imposed on
a simply connected subset of $\partial \oma$}.
Then, there holds
\be \label{eq_curl_Np_map_oma}
    \Crl (\ND_p(\Ta) \cap \Hca) = \{\tv_h \in \RT_p(\Ta) \cap \Hdva; \, \Dv \tv_h = 0\},
\ee
so that the minimization set in~\eqref{eq_bha}
is nonempty when
\be \label{eq_bha_min_set} \bs
& \frh^\ver + \RTproj{p} (\btha - \bdla) \in \RT_p(\Ta) \cap \Hdva \\
& \quad \text{ and } \quad \Dv (\frh^\ver + \RTproj{p} (\btha - \bdla))=0.
\es \ee

The first requirement in~\eqref{eq_bha_min_set} is immediate from~\eqref{eq_sigma_a}, \eqref{eq_tha}, \eqref{eq_bdla_div}, and~\eqref{eq_RT_proj}; please note the crucial use of the canonical elementwise Raviart--Thomas \cred{projector} $\RTproj{p}$ which
reduces the order $(p+1)$ down to $p$, while preserving the homogeneous Neumann boundary condition. We also stress that we crucially only employ $\RTproj{p}$ to the piecewise polynomials $\btha$ and $\bdla$. As for the divergence constraint in the second requirement in~\eqref{eq_bha_min_set}, we have from~\eqref{eq_sigma_a}
\[
    \Dv \frh^\ver = \Pi_h^p(\psia \Dv (\Crl \tv) + \Gr \psia \scp (\Crl \tv)) = \Pi_h^p(\Gr \psia \scp (\Crl \tv)).
\]
We have, in turn, from the commuting property~\eqref{eq_com_prop_div} and~\eqref{eq_tha}
\[
    \Dv \RTproj{p}( \btha) = \Pi_h^p(\Dv \btha) = \Pi_h^p(\Pi_h^{p+1}(- \Gr \psia \scp (\Crl \tv))) = \Pi_h^p(- \Gr \psia \scp (\Crl \tv)),
\]
whereas $\bdla$, and consequently $\RTproj{p}(\bdla)$, are divergence-free from~\eqref{eq_bdla_div} \cred{and~\eqref{eq_com_prop_div}}.

Finally, \eqref{eq_bh_prop} is an immediate consequence of $\bha \in \ND_p(\Ta) \cap \Hca$ and of the definition~\eqref{eq_bh}.
\ep

\section{Proof of Theorem~\ref{thm_Phc}} \label{sec_Phc_proof}

We prove the claims of Theorem~\ref{thm_Phc} separately. Let the assumptions of Theorem~\ref{thm_Phc} be satisfied.

\bl[Commuting property~\eqref{eq_com}] The commuting property~\eqref{eq_com} holds true. \el

\bp
We have
\ban
    \Crl \Phc (\tv) & \reff{eq_bh}= \suma \Crl \bha \reff{eq_bha}=\suma \big\{\frh^\ver + \RTproj{p} (\btha - \bdla)\big\} \\
    & = \suma \frh^\ver + \RTproj{p} \Bigg(\suma (\btha - \bdla) \Bigg)   \\ & \refff{eq_td}{eq_bdla_sum}= \suma \frh^\ver + \RTproj{p} (\bdl_h - \bdl_h) = \suma \frh^\ver \reff{eq_sigma}= \Phd (\Crl \tv),
\ean
where we also have crucially used the linearity of the canonical elementwise Raviart--Thomas \cred{projector} $\RTproj{p}$ from~\eqref{eq_RT_proj} in the third step.

\ep

\bl [Approximation property~\eqref{eq_approx}] \label{lem_approx} The approximation property~\eqref{eq_approx} holds true. \el

\bp Let a mesh element $\elm \in \Th$ be fixed. We proceed in two steps.

{\em Step~1}. We first bound the term $\frac{h_\elm}{p+1}\norm[\big]{\Crl\big(\tv-\Phc (\tv)\big)}_\elm$.
From the commuting property~\eqref{eq_com}, we have $\Crl \Phc (\tv) = \Phd (\Crl \tv)$, so that
\be \label{eq_est_curl} \bs
    \frac{h_\elm}{p+1}\norm[\big]{\Crl\big(\tv-\Phc (\tv)\big)}_\elm & = \frac{h_\elm}{p+1} \norm[\big]{\Crl\tv-\Phd (\Crl \tv)}_\elm \\
     & \ls \Bigg\{\sum_{\elmtt\in\TK} \bigg(\frac{h_{\elmtt}}{p+1}\norm[\big]{\Crl\tv - \RToproj{p}(\Crl\tv)}_\elmtt\bigg)^2 \Bigg\}^\ft
\es \ee
from~\cite[Theorem~3.2, bound~(3.6)]{Ern_Gud_Sme_Voh_loc_glob_div_22}, also using that $\Dv(\Crl\tv)=0$, the shape-regularity of the mesh yielding $h_{\elm} \approx h_{\elmtt}$, and the characterization~\eqref{eq_RToproj_eq}.

{\em Step~2}. We bound $\norm[\big]{\tv - \Phc (\tv)}_\elm$. Consider for this purpose $\biot_h$ defined in~\eqref{eq_bioth}. Using the linearity and projection property of the elementwise canonical N\'ed\'elec \cred{projector} $\NDproj{p}$ from~\eqref{eq_ND_proj}, the partition of unity~\eqref{eq_PU}, and definition~\eqref{eq_bh}
\ban
    \norm[\big]{\tv - \Phc (\tv)}_\elm & = \norm[\bigg]{\tv - \biot_h + \sum_{\ver \in \VK} \big\{ \NDproj{p} (\psia \biot_h) - \bha \big\}}_\elm \\
    & \leq \norm{\tv - \biot_h}_\elm + \sum_{\ver \in \VK} \norm[\big]{\NDproj{p} (\psia \biot_h) - \bha}_\oma.
\ean
We now bound the two terms separately.

{\em Step~2a}. We bound $\norm{\tv - \biot_h}_\elm$. Using the second inequality in~\eqref{eq_loc_const_eq} of Lemma~\ref{lem_loc_const_eq} below, we have
\be \label{eq_v_iota} \bs
    \norm{\tv - \biot_h}_\elm & =  \min_{\substack{\tv_h \in \ND_p(\elm)
\\
\Crl \tv_h = \btau_h
}}
\norm{\tv - \tv_h}_\elm \\
& \ls \Big(\min_{\tv_h \in \ND_p(\elm)}\norm{\tv - \tv_h}_\elm + \frac{h_{\elm}}{p+1}\norm[\big]{\Crl\tv - \RToproj{p}(\Crl\tv)}_\elm \Big),
\es \ee
for a constant actually only depending on the shape-regularity parameter $\kappa_{\Th}$.


{\em Step~2b}. Let a vertex $\ver \in \VK$. Then
Theorem~3.3 and Corollary~4.3 of~\cite{Chaum_Voh_p_rob_3D_H_curl_23} give
\be\label{eq_stab_patch}
\min_{\substack{
\tv_h \in \ND_p(\Ta) \cap \Hca
\\
\Crl \tv_h = \frh^\ver + \RTproj{p} (\btha - \bdla)
}}
\norm[\big]{\NDproj{p} (\psia \biot_h) - \tv_h}_\oma
\ls \min_{\substack{
\bvf \in \Hca
\\
\Crl \bvf = \frh^\ver + \RTproj{p} (\btha - \bdla)
}}
\norm[\big]{\NDproj{p} (\psia \biot_h) - \bvf}_\oma,
\ee
where again the constant beyond $\ls$ only depends on the shape-regularity parameter $\kappa_{\Th}$. This result builds on \cite{Demk_Gop_Sch_ext_II_09, Cost_McInt_Bog_Poinc_10} and~\cite{Chaum_Ern_Voh_curl_elm_20} and extends~\cite{Brae_Pill_Sch_p_rob_09, Ern_Voh_p_rob_3D_20, Chaum_Ern_Voh_Maxw_22} to vertex patches in $\HC$.

By the characterization~\eqref{eq_bha} for the above left-hand side and by a primal--dual equivalence for the above right-hand side  as in,
\eg, \cite[Lemma~5.5]{Chaum_Ern_Voh_Maxw_22}, we obtain
\be \label{eq_stab_est} \bs
    \norm[\big]{\NDproj{p} (\psia \biot_h) - \bha}_\oma \ls {} &
    \sup_{\substack{\bvf \in \Hcdual{\oma}\\ \norm{\Crl \bvf}_\oma = 1}}
    \big\{(\frh^\ver + \RTproj{p} (\btha - \bdla), \bvf)_\oma \\
    & - (\NDproj{p} (\psia \biot_h), \Crl \bvf)_\oma\big\}.
\es \ee
Fix now $\bvf \in \Hcdual{\oma}$ with $\norm{\Crl \bvf}_\oma = 1$. \cred{From Lemma~\ref{lem_broken_reg} in Appendix~\ref{sec_broken_reg} below, we know that there exists $\bpsi \in \Hcdual{\oma}$ with
\[
    \bpsi\vert_\elmtt \in \tHoi{\elmtt} \qquad \forall \elmtt \in \Ta
\]
such that
\bse \label{eq_broken_reg} \be
    \Crl \bpsi = \Crl \bvf \label{eq_curl_eq}
\ee
and
\ba
    \norm{\bpsi}_\oma & \ls h_\oma \norm{\Crl \bvf}_\oma = h_\oma, \label{eq_est_PF}\\
    \Bigg\{\sum_{\elmtt \in \Ta} \norm{\Gr \bpsi}_\elmtt^2 \Bigg\}^\ft & \ls \norm{\Crl \bvf}_\oma = 1. \label{eq_est_lift}
\ea \ese
}
In the second term on the right-hand side of~\eqref{eq_stab_est}, $\bvf$ can be exchanged for
$\bpsi$ \cred{using~\eqref{eq_curl_eq}}. This is also true for the first one: from~\eqref{eq_bha_min_set}
and~\eqref{eq_curl_Np_map_oma}, there exists $\ttw_h \in \ND_p(\Ta) \cap \Hca$ such that
$\Crl \ttw_h = \frh^\ver + \RTproj{p} (\btha - \bdla)$, so that
\be \label{eq_vol} \bs
(\frh^\ver + \RTproj{p} (\btha - \bdla), \bvf)_\oma
& =
(\Crl \ttw_h, \bvf)_\oma
=
(\ttw_h, \Crl \bvf)_\oma
\\
& =
(\ttw_h, \Crl \bpsi)_\oma
=
(\frh^\ver + \RTproj{p} (\btha - \bdla), \bpsi)_\oma.
\es \ee

\cred{Writing as a sum over the elements staring the vertex $\ver$ and} integrating \cred{elementwise} by parts, the second term on the right-hand side of~\eqref{eq_stab_est} becomes
\be \label{eq_IPP} \bs
    - (\NDproj{p} (\psia \biot_h), \Crl \bpsi)_\oma = {} & - \sum_{\elmtt \in \Ta} (\NDproj{p} (\psia \biot_h), \Crl \bpsi)_\elmtt \\
    = {} & - \sum_{\elmtt \in \Ta} (\Crl (\NDproj{p} (\psia \biot_h)), \bpsi)_\elmtt \\
    & {} + \sum_{\sd \in \Faint} \<\jump{\NDproj{p} (\psia \biot_h)\vp\tn_\sd},\bpsi\>_\sd,
\es \ee
where $\Faint$ denotes the interior faces of $\Ta$ plus, for boundary patches, the faces in $\overline \GN$; the other faces do not appear because of the presence of the hat function $\psia$ and/or due to the zero \cred{tangential} trace \cred{of $\bpsi \in \Hcdual{\oma}$} on $\gD$.

{\em Step~2bi}. We address the jump term in~\eqref{eq_IPP}.
The trace inequality $\norm{\bpsi}_\sd^2 \ls \norm{\Gr \bpsi}_\elmtt \norm{\bpsi}_\elmtt + h_\elmtt^{-1}\norm{\bpsi}_\elmtt^2$ for any $\elmtt \in \Th$ sharing $\sd$ \cred{and} the Cauchy--Schwarz \cred{and Young} inequalities lead to
\[
    \sum_{\sd \in \Faint} \norm{\bpsi}_\sd^2 \cred{\ls \sum_{\elmtt \in \Ta} \big( h_\elmtt \norm{\Gr \bpsi}_\elmtt^2 + h_\elmtt^{-1}\norm{\bpsi}_\elmtt^2 \big) \refff{eq_est_PF}{eq_est_lift}\ls h_\oma}.
\]
The stability of the elementwise N\'ed\'elec \cred{projector} and the fact that $\norm{\psia}_{\infty,\oma} = 1$ give, for $\sd \in \Faint$,
\[
    \norm{\jump{\NDproj{p} (\psia \biot_h)\vp\tn_\sd}}_\sd \ls \norm{\jump{\biot_h \vp\tn_\sd}}_\sd,
\]
where henceforth $\ls$ adds a dependence on the polynomial degree $p$.
The above right-hand side is now in a form of a standard residual-based
\cred{a posteriori error} estimator. Introduce the face bubble function
$\psiF$ as a product of the three \cred{hat functions} $\psia$ of the
vertices of $\sd$; $\psiF$ is supported on the \cred{face patch subdomain}
$\omF$, which corresponds to the (two or one) tetrahedra sharing $\sd$, and is
\cred{zero on $\pt \omF \setminus \sd$}. Denote $\bzeta_\sd \eq \jump{\biot_h \vp\tn_\sd}$
\cred{on the face $\sd$} together with its extension to $\omF$
\cred{by constants in the direction vertex opposite to the face $\sd$ -- barycenter of $\sd$}
and note that
\[
    (\Crl \tv, \bzeta_\sd \psiF)_\omF - (\tv, \Crl (\bzeta_\sd \psiF))_\omF = 0.
\]
\cred{Indeed, this follows} since $\bzeta_\sd \psiF \in \tHoi{\omF}$, $\bzeta_\sd \psiF$ is zero on $\pt \omF \setminus \sd$, and $\tv \in \HcN$. Then, using standard scaling arguments as in, e.g., \cite[page~172]{Beck_Hipt_Hopp_Wohl_a_post_Maxw_00}, leads to
\ban
    \norm{\jump{\biot_h \vp\tn_\sd}}_\sd^2 & = \norm{\bzeta_\sd}_\sd^2 \ls \<\jump{\biot_h \vp\tn_\sd}, \bzeta_\sd \psiF\>_\sd \\
    & = (\tv - \biot_h, \Crl (\bzeta_\sd \psiF))_\omF - (\Crl (\tv - \biot_h), \bzeta_\sd \psiF)_\omF\\
    & \leq \norm{\tv - \biot_h}_\omF \norm{\Crl (\bzeta_\sd \psiF)}_\omF + (\btau_h - \Crl \tv, \bzeta_\sd \psiF)_\omF,
\ean
where we have also used that $\Crl \biot_h = \btau_h$ from~\eqref{eq_bioth}. The second term above is treated by Lemma~\ref{lem_hp_d_est} below with $\bpsi = \bzeta_\sd \psiF$ on each element \cred{$\elmtt$} in the face patch subdomain $\omF$. This yields
\[
    (\btau_h - \Crl \tv, \bzeta_\sd \psiF)_\omF \ls \frac{h_\omF}{p+1}\norm[\big]{\Crl\tv - \RToproj{p}(\Crl\tv)}_\omF \norm{\Gr (\bzeta_\sd \psiF)}_\omF.
\]
By the inverse inequality,
\[
    \norm{\Crl (\bzeta_\sd \psiF)}_\omF \ls h_\omF^{-1} \norm{\bzeta_\sd \psiF}_\omF, \quad
    \norm{\Gr (\bzeta_\sd \psiF)}_\omF \ls h_\omF^{-1} \norm{\bzeta_\sd \psiF}_\omF,
\]
whereas the constant extension of $\bzeta_\sd$ from $\sd$ to $\omF$ and $\norm{\psiF}_{\infty,\omF} \leq 1$ give
\[
    \norm{\bzeta_\sd \psiF}_\omF \ls h_\omF^\ft \norm{\bzeta_\sd}_\sd.
\]
Thus, altogether,
\[
    h_\omF^\ft \norm{\bzeta_\sd}_\sd \ls \norm{\tv - \biot_h}_\omF + \frac{h_\omF}{p+1}\norm[\big]{\Crl\tv - \RToproj{p}(\Crl\tv)}_\omF.
\]
This yields the desired components in~\eqref{eq_approx}, using the mesh shape regularity and the definition~\eqref{eq_bioth} of $\biot_h$ together with the constrained--unconstrained equivalence result of Lemma~\ref{lem_loc_const_eq} below, which allows to pass from $\norm{\tv - \biot_h}_\elmtt$ to $\min_{\tv_h \in \ND_p(\elmtt)}\norm{\tv - \tv_h}_\elmtt$ plus the oscillation terms on each simplex $\elmtt$\cred{, as in~\eqref{eq_v_iota}}.

{\em Step~2bii}. We now address the volume term in~\eqref{eq_IPP} together with~\eqref{eq_vol}. On each $\elmtt \in \Ta$, using the commuting property of the canonical \cred{projector}~\eqref{eq_com_prop_curl},
\be \label{eq_ipN}
    \Crl (\NDproj{p} (\psia \biot_h)) = \RTproj{p}(\Crl(\psia \biot_h)) = \RTproj{p}(\Gr \psia \vp \biot_h + \psia \Crl \biot_h).
\ee

Recall~\eqref{eq_RToproj_eq} and~\eqref{eq_sigma_a}, where $\ttw = \Crl \tv$, and Remark~\ref{rem_contr_mod}. Thus, from~\cite[Lemma~4.6]{Ern_Gud_Sme_Voh_loc_glob_div_22}, we have the bound
\be \label{eq_sigma_a_bound}
    \norm{\frh^\ver - \RTproj{p}(\psia \btau_h)}_\oma \ls \Bigg\{\sum_{\elmtt\in\Ta} \norm{\Crl\tv - \RToproj{p}(\Crl\tv)}_\elmtt^2\Bigg\}^\ft,
\ee
where we note that the data oscillation terms disappears since $\Dv(\Crl \tv)=0$.
Thus, using~\eqref{eq_ipN} and regrouping the terms in~\eqref{eq_IPP} and~\eqref{eq_vol}, since elementwise $\Crl \biot_h = \btau_h$ from~\eqref{eq_bioth}, and using the Cauchy--Schwarz inequality together with~\eqref{eq_est_PF} and the mesh shape regularity,
\ban
    (\frh^\ver - \RTproj{p}(\psia \Crl \biot_h), \bpsi)_\oma & = (\frh^\ver - \RTproj{p}(\psia \btau_h), \bpsi)_\oma \\
    & \leq \norm{\frh^\ver - \RTproj{p}(\psia \btau_h)}_\oma \norm{\bpsi}_\oma \\
    & \ls \Bigg\{\sum_{\elmtt\in\Ta} \norm{\Crl\tv - \RToproj{p}(\Crl\tv)}_\elmtt^2\Bigg\}^\ft h_\oma \\
    & \ls \Bigg\{\sum_{\elmtt\in\Ta} \bigg(\frac{h_{\elmtt}}{p+1}\norm{\Crl\tv - \RToproj{p}(\Crl\tv)}_\elmtt\bigg)^2\Bigg\}^\ft,
\ean
where we have also used \cred{the non-$p$-robust estimate} $1 \ls 1/(p+1)$.

We now start to estimate the remaining term. By the Cauchy--Schwarz inequality together with~\eqref{eq_est_PF}, we have
\be \label{eq_last} \bs
    (\RTproj{p} (\btha - \bdla - \Gr \psia \vp \biot_h), \bpsi)_\oma & \ls [\norm{\btha - \Gr \psia \vp \biot_h}_\oma + \norm{\bdla}_\oma] \norm{\bpsi}_\oma \\
    & \ls h_\oma [\norm{\btha - \Gr \psia \vp \biot_h}_\oma + \norm{\bdla}_\oma],
\es \ee
where we have also employed the stability of the elementwise Raviart--Thomas \cred{projector}.
We now employ definition~\eqref{eq_tha} and, crucially, \cite[Theorem~A.2]{Chaum_Voh_Maxwell_equil_23} (\cf\ also Lemma~7.3 therein for the rewriting of the data oscillation term). We have already verified in the proof of Lemma~\ref{lem_tha} that our setting fits the assumptions of~\cite[Theorem~A.2]{Chaum_Voh_Maxwell_equil_23}. We additionally note that $\Gr \psia \vp \tv$ lies in $\Hdva$ with $\Dv (\Gr \psia \vp \tv) = - \Gr \psia \scp (\Crl \tv)$. This leads to
\be \label{eq_tha_eff} \bs
    & \norm{\btha - \Gr \psia \vp \biot_h}_\oma \\
    \ls {} & \norm{\Gr \psia \vp (\tv - \biot_h)}_\oma + h_\oma^{-1} \Bigg\{\sum_{\elmtt \in \Ta} \Big(\frac{h_\elmtt}{\pi} \norm{\Crl \tv - \bPi_h^{p+1}(\Crl \tv)}_\elmtt\Big)^2\Bigg\}^\ft \\
    \ls {} & h_\oma^{-1} \norm{\tv - \biot_h}_\oma + h_\oma^{-1}\Bigg\{\sum_{\elmtt\in\Ta} \bigg(\frac{h_{\elmtt}}{p+1}\norm{\Crl\tv - \RToproj{p}(\Crl\tv)}_\elmtt\bigg)^2\Bigg\}^\ft.
\es \ee
Thus, definition~\eqref{eq_bioth} of $\biot_h$ together with the constrained--unconstrained equivalence result of Lemma~\ref{lem_loc_const_eq} below \cred{as in~\eqref{eq_v_iota}} yields the desired bound for the first term in~\eqref{eq_last}.

Finally, we bound $\norm{\bdla}_\oma$. The construction of $\bdl_h$ and $\bdla$ from steps~\ref{dfc_3} and~\ref{dfc_4} of Definition~\ref{def_Phc} fits the framework of~\cite[Theorem~B.1]{Chaum_Voh_Maxwell_equil_23} with $q=q'=p+1$, as we have already verified in the proof of Lemma~\ref{lem_delta}. Thus, using~\cite[estimate~(B.6b)]{Chaum_Voh_Maxwell_equil_23}, \eqref{eq_td}, and~\eqref{eq_PU}, we have, for any $\elmtt\in\Ta$,
\[
    \norm{\bdla}_\elmtt \ls \norm{\bdl_h}_\elmtt = \norm[\Bigg]{\sum_{\vertt\in\VKtt}(\bthb-\Gr \psib \vp \biot_h)}_\elmtt,
\]
so that the bound on $\norm{\bdla}_\oma$ follows from that on $\norm{\bthb - \Gr \psib \vp \biot_h}_\omb$ (this final estimate extends the bound~\eqref{eq_approx} from neighbors of $\elm$ to the neighbors of neighbors of $\elm$).

\ep

\bl [Projection property~\eqref{eq_proj}] The projection property~\eqref{eq_proj} holds true.\el

\bp When $\tv \in \ND_p(\Th) \cap \HcN$ is a N\'ed\'elec piecewise polynomial, the right-hand side in the approximation property~\eqref{eq_approx} is zero, so that immediately $\Phc (\tv) = \tv$ follows. \ep

\bl[Stability properties~\eqref{eq_stab} and~\eqref{eq_stab_curl_norm}] The stability properties~\eqref{eq_stab} and~\eqref{eq_stab_curl_norm} hold true.\el

\bp The triangle inequality gives
\[
    \norm[\big]{\Phc (\tv)}_\elm \leq \norm[\big]{\tv - \Phc (\tv)}_\elm + \norm{\tv}_\elm.
\]
Thus~\eqref{eq_stab} follows from~\eqref{eq_approx} and the orthogonal projection stability
\be \label{eq_ND_stab}
    \min_{\tv_h \in \ND_p(\elmtt)}\norm{\tv - \tv_h}_\elmtt \leq \norm{\tv}_\elmtt.
\ee

As for~\eqref{eq_stab_curl_norm}, we only need to treat the second term. Indeed~\eqref{eq_stab} bounds the first term as $h_\elmtt/(p+1) \leq h_\Om$ \cred{and since,} as in~\eqref{eq_ND_stab},
\be \label{eq_RT_stab}
    \norm[\big]{\Crl\tv - \RToproj{p}(\Crl\tv)}_\elmtt \reff{eq_RToproj_eq}= \min_{\ttw_h \in \RT_p(\elmtt)} \norm{\Crl\tv - \ttw_h}_\elm  \leq \norm{\Crl\tv}_\elmtt.
\ee
The triangle inequality gives
\[
    h_\Om \norm[\big]{\Crl \Phc (\tv)}_\elm \leq h_\Om \norm[\big]{\Crl (\tv - \Phc (\tv))}_\elm + h_\Om \norm[\big]{\Crl \tv}_\elm,
\]
and the first term above is estimated as in~\eqref{eq_est_curl}, with the weight $h_\Om$ in place of $h_\elm/(p+1)$, \cred{and employing~\eqref{eq_RT_stab}.} \ep

\section{Proof of Theorem~\ref{thm_loc_glob} and of~\eqref{eq_LG_constr_glob}} \label{sec_loc_glob_proof}

The proof of properties~\eqref{eq_LG_constr} and~\eqref{eq_LG_unconstr} follows straightforwardly from Theorem~\ref{thm_Phc}, whereas the proof~\eqref{eq_LG_constr_glob} will turn out slightly more involved.

\bp[Proof of Theorem~\ref{thm_loc_glob}.] {\em Step~1}. We first show~\eqref{eq_LG_constr}.
Since the second terms are identically present on both sides of the equivalence (\cf~\cite[Remark~3.4]{Ern_Gud_Sme_Voh_loc_glob_div_22}), we only need to consider the first ones. As the inequality trivially holds with the $\geq$ sign, we only need to bound the first term of the left hand side of~\eqref{eq_LG_constr} from above.
\cred{Employing the projector $\Phc$ from Definition~\ref{def_Phc} and the commuting property~}\eqref{eq_com}, we immediately have
\be \label{eq_b1}
    \min_{\substack{\tv_h \in \ND_p(\Th) \cap \HcN\\\Crl\tv_h = \Phd (\Crl \tv)}}\norm{\tv - \tv_h}^2 \leq \norm[\big]{\tv - \Phc (\tv)}^2,
\ee
so that the result follows from the local approximation property~\eqref{eq_approx} by summing over all mesh elements $\elm \in \Th$ and invoking the mesh shape-regularity.

{\em Step~2}. We now show~\eqref{eq_LG_unconstr}. This again trivially holds with $\geq$, since
\be \label{eq_curl_Np_map}
\Crl (\ND_p(\Th) \cap \HcN) \subset \{\tv_h \in \RT_p(\Th) \cap \HdvN; \, \Dv \tv_h = 0\};
\ee
thus both terms on the right-hand side of~\eqref{eq_LG_unconstr} are the elementwise versions of those on the left. For the converse estimate, we bound the left-hand side of~\eqref{eq_LG_unconstr} by
\[
    \norm{\tv - \Phc (\tv)}^2 + \sum_{\elm \in \Th}\bigg(\frac{h_{\elm}}{p+1}\norm{\Crl(\tv-\Phc (\tv))}_\elm\bigg)^2
\]
and again invoke~\eqref{eq_approx}. \ep

\bp[Proof of~\eqref{eq_LG_constr_glob}.]

{\em Step~1}. We recall that here,
$p \geq 1$, $\Th$ is quasi-uniform, and $\Omega$ is convex with either $\GD = \emptyset$ or $\GN = \emptyset$.
The second term on the left-hand side of~\eqref{eq_LG_constr_glob} is again the same as on the right-hand side, so we only need to treat the first one. Let
\be \label{eq_zh_def}
    \bzeta_h \eq \arg \min_{\substack{\tv_h \in \ND_p(\Th) \cap \HcN\\\Crl\tv_h = \bPi_h^{p,\dv}(\Crl \tv) - \Phd (\Crl \tv)}}\norm{\tv_h}^2.
\ee
Note that $\bzeta_h$ is only nonzero when the datum $\tv \in \HcN$ has a non-polynomial curl, $\Crl \tv \not \in [\PP_p(\Th)]^3$ (or\cred{, equivalently,} $\Crl \tv \not \in \RT_p(\Th) \cap \HdvN$).
Since $\Phc(\tv) + \bzeta_h$ lies in $\ND_p(\Th) \cap \HcN$ and satisfies $\Crl(\Phc(\tv) + \bzeta_h) = \bPi_h^{p,\dv}(\Crl \tv)$, it is clear that
\[
    \min_{\substack{\tv_h \in \ND_p(\Th) \cap \HcN\\\Crl\tv_h = \bPi_h^{p,\dv}(\Crl \tv)}}\norm{\tv - \tv_h} \leq \norm{\tv - \Phc(\tv)} + \norm{\bzeta_h}.
\]
The first term on the right-hand side above is as in~\eqref{eq_b1}, so we only treat the second one.

Let
\[
    \bzeta \eq \arg \min_{\substack{\tv \in \HcN\\\Crl\tv = \bPi_h^{p,\dv}(\Crl \tv) - \Phd (\Crl \tv)}}\norm{\tv}^2;
\]
$\bzeta$ is again only nonzero when the datum $\tv \in \HcN$ has a non-polynomial curl.
Since the curl constraint in the above minimization problem is given by a Raviart--Thomas piecewise polynomial, we infer from the definition of $\bzeta_h$, the commuting property~\eqref{eq_com}, and the stability~\eqref{eq_stab} that
\[
    \norm{\bzeta_h} \leq \norm{\Phc(\bzeta)} \ls \norm{\bzeta},
\]
where the constant hidden in $\ls$ only depends on the shape-regularity parameter $\kappa_{\Th}$ of the mesh $\Th$ and the polynomial degree $p$. By a primal--dual equivalence, as in,
\eg, \cite[Lemma~5.5]{Chaum_Ern_Voh_Maxw_22}, we have
\[
    \norm{\bzeta} = \sup_{\substack{\bvf \in \HcD\\ \norm{\Crl \bvf} = 1}}
    \big\{\big(\bPi_h^{p,\dv}(\Crl \tv) - \Phd (\Crl \tv), \bvf\big)\big\}.
\]
Let now
\be \label{eq_bvf_Om}
    \bvf \in \HcD \text{ with } \norm{\Crl \bvf} = 1
\ee
be fixed. Define $q \in \HoD$ by
\be \label{eq_q_Om}
    (\Gr q, \Gr v) = (\cred{\bvf}, \Gr v) \qquad \forall v \in \HoD.
\ee
Defining $\bchi \eq \cred{\bvf} - \Gr q$, we have, recalling the notation from Section~\ref{sec_cont_sp_BC},
\be \label{eq_bchi_reg_Om}
    \bchi \in \HcD \cap \HdvN \text{ with } \Dv \bchi = 0,
\ee
together with
\be \label{eq_crl_bxi}
    \Crl \bchi = \Crl \bvf.
\ee
Moreover, since $\Om$ is supposed convex here and either $\GD = \emptyset$ or $\GN = \emptyset$,
\cred{we can employ~\cite[Theorems~3.7 or~3.9]{Gir_Rav_NS_86}. Thus there also holds $\bchi \in \tHo$ with}
\be \label{eq_gr_chi_Om}
    \norm{\Gr \bchi} \ls \norm{\Crl \bchi} \reff{eq_crl_bxi}= \norm{\Crl \bvf} \reff{eq_bvf_Om}=1.
\ee
\cred{[--]}.
Since \cred{from~\eqref{eq_zh_def}} $\Crl\bzeta_h  = \bPi_h^{p,\dv}(\Crl \tv) - \Phd (\Crl \tv)$, we can exchange $\bvf$ for $\cred{\bchi}$, so that
\ban
    \big(\bPi_h^{p,\dv}(\Crl \tv) - \Phd (\Crl \tv), \bvf\big)
    = {} & \big(\bPi_h^{p,\dv}(\Crl \tv) - \Phd (\Crl \tv), \cred{\bchi}\big) \\
    = {} & \big(\bPi_h^{p,\dv}(\Crl \tv) - \Crl \tv, \cred{\bchi}\big) \\
    {} & + \big(\Crl \tv - \Phd (\Crl \tv), \cred{\bchi}\big).
\ean
We will treat the two arising terms separately.

{\em Step~2}. \cred{[--]}.
Denote henceforth $\ttw \eq \Crl \tv$. From the definition~\eqref{eq_RT_glob} of the projector $\bPi_h^{p,\dv}$, we can subtract any $\bchi_h \in \RT_p(\Th) \cap \HdvN$ with $\Dv \bchi_h = 0$; we choose $\bchi_h \eq \Phd (\bchi)$ following Definition~\ref{def_Phd}, yielding
\ban
    \big(\bPi_h^{p,\dv}(\ttw) - \ttw, \bchi\big) & = \big(\bPi_h^{p,\dv}(\ttw) - \ttw, \bchi - \Phd(\bchi)\big) \\
    & \leq \norm[\big]{\bPi_h^{p,\dv}(\ttw) - \ttw} \norm[\big]{\bchi - \Phd(\bchi)}.
\ean
Now the equivalence of global-best and local-best approximations from~\cite[Theorem~3.3]{Ern_Gud_Sme_Voh_loc_glob_div_22} gives, using that $\Dv \ttw = 0$ \cred{and recalling $\ttw = \Crl \tv$},
\[
    \norm[\big]{\bPi_h^{p,\dv}(\ttw) - \ttw} \ls \Bigg\{\sum_{\elm\in\Th} \norm{\Crl\tv - \RToproj{p}(\Crl\tv)}_\elm^2\Bigg\}^\ft,
\]
where again the constant of $\ls$ only depends on the shape-regularity parameter $\kappa_{\Th}$ of the mesh $\Th$ and the polynomial degree $p$. In particular, it seems impossible to estimate $\norm[\big]{\bPi_h^{p,\dv}(\ttw) - \ttw}_\elm$, locally on each $\elm \in \Th$. Similarly, the approximation properties of the projector $\Phd$ from equation~(3.6) of~\cite[Theorem~3.2]{Ern_Gud_Sme_Voh_loc_glob_div_22} \cred{and using $\bchi \in \tHo$} lead to
\ban
    \norm[\big]{\bchi - \Phd(\bchi)} & \ls \Bigg\{\sum_{\elm\in\Th} \min_{\ttw_h \in \RT_p(\elm)}
    \norm{\bchi - \ttw_h}_\elm^2 \Bigg\}^\ft \\
    & \ls \Bigg\{\sum_{\elm\in\Th} h_\elm^2 \norm{\Gr \bchi}_\elm^2 \Bigg\}^\ft \leq h \norm{\Gr \bchi} \cred{\reff{eq_gr_chi_Om} \ls h},
\ean
where, recall, $h \eq \max_{\elm \in \Th} h_\elm$.
Thus, since $\Th$ is supposed here quasi-uniform, so that $h \ls h_\elm$ for all $\elm \in \Th$, and using \cred{the non-$p$-robust estimate} $1 \ls 1/(p+1)$ as in the proof of Lemma~\ref{lem_approx}, we conclude
\[
    \big(\bPi_h^{p,\dv}(\Crl \tv) - \Crl \tv, \bchi\big) \ls \Bigg\{\sum_{\elm\in\Th} \bigg(\frac{h_{\elm}}{p+1}\norm{\Crl\tv - \RToproj{p}(\Crl\tv)}_\elm\bigg)^2\Bigg\}^\ft.
\]

{\em Step~3}. Denote $\ttw \eq \Crl \tv$ as above. Consider $\btau_h$ given by~\eqref{eq_tauh}\cred{, so that in particular $\Dv \btau_h = 0$}. Using the partition of unity~\eqref{eq_PU}, the linearity and the projection property of the elementwise Raviart--Thomas \cred{projector} $\RTproj{p}$ of~\eqref{eq_RT_proj}, we have
\[
    \btau_h = \suma (\psia \btau_h) = \RTproj{p}(\btau_h) = \suma \big(\RTproj{p}(\psia \btau_h)\big).
\]
Thus, also using~\eqref{eq_sigma},
\ban
    \big(\ttw - \Phd (\ttw), \cred{\bchi}\big) & = \suma \big(\psia \ttw - \psia \btau_h + \RTproj{p}(\psia \btau_h) - \frh^\ver, \cred{\bchi}\big)_\oma.
\ean
Consider now any $q_h \in \PP_1(\Ta) \cap \Hsa$ for a fixed \cred{vertex} $\ver \in \Vh$. Since $\ttw \in \HdvN$ with $\Dv \ttw = 0$ and by the constraint in~\eqref{eq_sigma_a}, we see, \cred{on the one hand,} employing the Green theorem,
\ban
    \big(\psia \ttw - \frh^\ver, \Gr q_h\big)_\oma & = - \big(\Dv(\psia \ttw - \frh^\ver), q_h\big)_\oma \\
    & = - \big(\Gr \psia \scp \ttw - \Pi_h^p (\Gr \psia \scp \ttw)), q_h\big)_\oma = 0.
\ean
On the other \cred{hand}, on each simplex in the patch $\elm \in \Ta$, the Green theorem, \eqref{eq_RT_proj}, and the commuting property~\eqref{eq_com_prop_div} of $\RTproj{p}$ give
\ban
    \big(\psia \btau_h - \RTproj{p}(\psia \btau_h), \Gr q_h\big)_\elm = {} & - \big(\Dv\big(\psia \btau_h - \RTproj{p}(\psia \btau_h)\big), q_h\big)_\elm \\
    & {} + \big\<\big(\psia \btau_h -  \RTproj{p}(\psia \btau_h)\big) \scp \tn_\elm, q_h\big\>_{\pt \elm},\\
    = {} & - \big(\Gr \psia \scp \btau_h - \Pi_h^p(\Gr \psia \scp \btau_h), q_h\big)_\elm = 0.
\ean
For the last two properties, we have also employed the assumption $p \geq 1$.

We now choose $q_h \in \PP_1(\Ta) \cap \Hsa$ such that $\Gr q_h = \cred{\overline {\bchi\lvert_\oma}}$, the componentwise mean value of $\cred{\bchi}$ on the patch subdomain $\oma$.
Thus, the Cauchy--Schwarz inequality and the Poincar\'e--Friedrichs inequality~\eqref{eq_Poinc_Fried} give
\ban
    & \big(\psia \ttw - \psia \btau_h + \RTproj{p}(\psia \btau_h) - \frh^\ver, \cred{\bchi}\big)_\oma \\
    = {} & \big(\psia \ttw - \psia \btau_h + \RTproj{p}(\psia \btau_h) - \frh^\ver, \cred{\bchi} - \cred{\overline {\bchi\lvert_\oma}}\big)_\oma \\
    \ls {} & h_\oma
    \big(\norm{\psia \ttw - \psia \btau_h}_\oma + \norm[\big]{\RTproj{p}(\psia \btau_h) - \frh^\ver}_\oma \big) \norm{\Gr \cred{\bchi}}_\oma,
\ean
\cred{where we have again used $\bchi \in \tHo$.} Now, from the definition~\eqref{eq_tauh} of $\btau_h$ as a local constrained minimizer, from the
constrained--unconstrained equivalence of~\cite[Lemma~A1]{Ern_Gud_Sme_Voh_loc_glob_div_22} (recall that $\Dv \ttw=0$) \cred{(\cf~\eqref{eq_proj_est} below)}, and redeveloping $\ttw = \Crl \tv$, we see
\[
    \norm{\psia \ttw - \psia \btau_h}_\oma \leq \norm{\ttw - \btau_h}_\oma \ls \Bigg\{\sum_{\elmtt\in\Ta} \norm{\Crl\tv - \RToproj{p}(\Crl\tv)}_\elmtt^2\Bigg\}^\ft.
\]
Similarly, it follows from~\eqref{eq_sigma_a_bound} that
\[
    \norm[\big]{\RTproj{p}(\psia \btau_h) - \frh^\ver}_\oma \ls \Bigg\{\sum_{\elmtt\in\Ta} \norm{\Crl\tv - \RToproj{p}(\Crl\tv)}_\elmtt^2\Bigg\}^\ft.
\]
Hence, the shape-regularity of the mesh yielding finite overlap between the patches, the Cauchy--Schwarz inequality, and using \cred{the non-$p$-robust estimate} $1 \ls 1/(p+1)$ as in the proof of Lemma~\ref{lem_approx}, we see
\ban
    \big(\Crl \tv - \Phd (\Crl \tv), \cred{\bchi}\big) & \refvs{}{\ls}{xxx} \Bigg\{\sum_{\elm\in\Th} \bigg(\frac{h_{\elm}}{p+1}\norm{\Crl\tv - \RToproj{p}(\Crl\tv)}_\elm\bigg)^2\Bigg\}^\ft \norm{\Gr \cred{\bchi}} \\
    & \cred{\refvs{\eqref{eq_gr_chi_Om}}{\ls}{xxx} \Bigg\{\sum_{\elm\in\Th} \bigg(\frac{h_{\elm}}{p+1}\norm{\Crl\tv - \RToproj{p}(\Crl\tv)}_\elm\bigg)^2\Bigg\}^\ft}.
\ean
This finishes the proof. \ep

\section{Proof of Theorem~\ref{thm_hp}} \label{sec_hp_proof}

We proceed in the spirit of~\cite[Section~5]{Ern_Gud_Sme_Voh_loc_glob_div_22}. The main point is to derive an estimate of the form of $\ls$ in~\eqref{eq_LG_unconstr} which is $p$-robust, \ie, where the hidden constant is independent of the polynomial degree $p$, for the price of lowering the approximation polynomial degree on the right from $p$ to $p-1$.
For this purpose, we only consider $p \geq 1$ and present alternatives of Definitions~\ref{def_Phd} and~\ref{def_Phc} where all instances of the elementwise Raviart--Thomas and N\'ed\'elec \cred{projector}s $\RTproj{p}$, $\RTproj{p+1}$, and $\NDproj{p}$ are removed. In addition, we simplify Definition~\ref{def_Phc}, as commuting will not be needed here.

We start with reworking Definition~\ref{def_Phd}:

\bd[Alternative of Definition~\ref{def_Phd}] \label{def_Phd_II} Let $\tv \in \HcN$ be given. Let $\ttw \eq \Crl \tv \in \HdvN$ with $\Dv \ttw = 0$\cred{, a simplicial mesh $\Th$ of $\Om$, and a} polynomial degree $p \geq 1$ \cred{be given}:

\begin{enumerate} [parsep=0.2pt, itemsep=0.2pt, topsep=1pt, partopsep=1pt, leftmargin=13pt]

\item \label{dfd_1_II} Define a broken Raviart--Thomas polynomial $\btau_h \in \RT_{p-1}(\Th)$, on each mesh element~\cred{$\elm$}, via
\be \label{eq_tauh_II}
\btau_h\vert_\elm \eq \arg \min_{\substack{
\ttw_h \in \RT_{p-1}(\elm)
\\
\Dv \ttw_h = 0
}}
\norm{\ttw - \ttw_h}_\elm \qquad \forall \elm \in \Th.
\ee

\item \label{dfd_2_II} Define a Raviart--Thomas polynomial $\frh^\ver \in \RT_{p}(\Ta) \cap \Hdva$, on each vertex patch~\cred{$\Ta$}, via
\be \label{eq_sigma_a_II}
\frh^\ver
\eq
\arg \min_{\substack{
\ttw_h \in \RT_{p}(\Ta) \cap \Hdva
\\
\Dv \ttw_h= \Pi_h^{p}(\Gr \psia \scp \ttw)
}}
\norm[\big]{\psia \btau_h - \ttw_h}_\oma \qquad \forall \ver \in \Vh.
\ee

\item \label{dfd_3_II} Extending $\frh^\ver$ by zero outside of the patch subdomain $\oma$, define $\frh \in \RT_{p}(\Th) \cap \HdvN$ via
\be \label{eq_sigma_II}
    \frh \eq \suma \frh^\ver.
\ee

\end{enumerate}
\ed

Definition~\ref{def_Phd_II} takes the same form as~\cite[Definition~5.2]{Ern_Gud_Sme_Voh_loc_glob_div_22}, up to the fact that the divergence constraint in~\eqref{eq_sigma_a_II} is modified in that we use $\Pi_h^{p}(\Gr \psia \scp \ttw)$ in place of $\Gr \psia \scp \btau_h$. Similarly to Remark~\ref{rem_contr_mod}, we now inspect the (beginning of) the proof of~\cite[Lemma~5.3]{Ern_Gud_Sme_Voh_loc_glob_div_22}.
Here a term $(g_\ver,\vf)_\oma = (\Pi_{\mathcal{T}}^p(\psia \Dv \tv) + \Gr \psia \scp \btau_{\mathcal{T}})_\oma$ appears (notation from the proof of~\cite[Lemma~5.3]{Ern_Gud_Sme_Voh_loc_glob_div_22}), which has to be replaced by $(\Pi_{\mathcal{T}}^p(\psia \Dv \tv + \Gr \psia \scp \tv), \vf)_\oma$. This makes appear the same supplementary term and its bound as in~\eqref{eq_sup_term}. Thus, \cite[Lemma~5.3]{Ern_Gud_Sme_Voh_loc_glob_div_22}, recalling~\eqref{eq_RToproj_eq} and $\Dv \ttw = \Dv (\Crl \tv) = 0$ (current notation), so that there are no oscillation terms, lead to
\be \label{eq_div_p_1}
    \norm[\big]{\psia \btau_h - \frh^\ver}_\oma \ls \Bigg\{\sum_{\elm\in\Ta} \norm{\Crl\tv - \RToproj{p-1}(\Crl\tv)}_\elm^2\Bigg\}^\ft,
\ee
where crucially the hidden constant only depends on the shape-regularity parameter $\kappa_{\Th}$.

We now rework Definition~\ref{def_Phc}.

\bd[Alternative of Definition~\ref{def_Phc}] \label{def_Phc_II} Let $\tv \in \HcN$\cred{, a simplicial mesh $\Th$ of $\Om$, and a polynomial degree} $p \geq 1$ be given:

\begin{enumerate} [parsep=0.2pt, itemsep=0.2pt, topsep=1pt, partopsep=1pt, leftmargin=13pt]

\item \label{dfc_0_II} For $\ttw \eq \Crl \tv$ yielding $\ttw \in \HdvN$ with $\Dv \ttw =0$,
define $\btau_h \in \RT_{p-1}(\Th)$ by~\eqref{eq_tauh_II} and $\frh^\ver \in \RT_p(\Ta) \cap \Hdva$ by~\eqref{eq_sigma_a_II} from Definition~\ref{def_Phd_II}.

\item \label{dfc_1_II} Define a broken N\'ed\'elec polynomial $\biot_h \in \ND_{p-1}(\Th)$, on each mesh element~\cred{$\elm$}, via
\be \label{eq_bioth_II}
\biot_h\vert_\elm \eq \arg \min_{\substack{
\tv_h \in \ND_{p-1}(\elm)
\\
\Crl \tv_h = \btau_h
}}
\norm{\tv - \tv_h}_\elm \qquad \forall \elm \in \Th.
\ee

\item \label{dfc_2_II} Define a Raviart--Thomas polynomial $\btha \in \RT_{p}(\Ta) \cap \Hdva$, on each vertex patch~\cred{$\Ta$}, via
\be \label{eq_tha_II}
\btha \eq \arg \hspace*{-1cm}
\min_{\substack{
\tv_h \in \RT_{p}(\Ta) \cap \Hdva
\\
\Dv \tv_h= \Pi_h^{p}(- \Gr \psia \scp (\Crl \tv))
}} \hspace*{-1cm} \norm{\Gr \psia \vp \biot_h- \tv_h}_\oma \qquad \forall \ver \in \Vh.
\ee

\item \label{dfc_5_II} Define a N\'ed\'elec polynomial $\bha \in \ND_p(\Ta) \cap \Hca$, on each vertex patch, via
\be \label{eq_bha_II}
\bha
\eq
\arg \min_{\substack{
\tv_h \in \ND_p(\Ta) \cap \Hca
\\
\Crl \tv_h = \frh^\ver + \btha
}}
\norm[\big]{\psia \biot_h - \tv_h}_\oma \qquad \forall \ver \in \Vh.
\ee

\item \label{dfc_6_II} Extending $\bha$ by zero outside of the patch subdomain $\oma$, define
\be \label{eq_bh_II}
    \bh_h \eq \suma \bha \in \ND_p(\Th) \cap \HcN.
\ee

\end{enumerate}

\ed

Compared to Definition~\ref{def_Phc}, Definition~\ref{def_Phc_II} avoids the construction of $\bdl_h$ on Step~\ref{dfc_3} and of the correction terms $\bdla$ on Step~\ref{dfc_4}. As in~\eqref{eq_td}, we can still set
\be \label{eq_td_II}
    \bdl_h \eq \suma \btha \in \RT_p(\Th) \cap \HdvN \quad \text{ with } \quad \Dv \bdl_h = 0,
\ee
but the commuting property $\Crl \bh_h = \frh$ from~\eqref{eq_com} is here lost in that from~\eqref{eq_bh_II}, \eqref{eq_bha_II}, \eqref{eq_sigma_II}, \cred{and~\eqref{eq_td_II}}

\be\label{eq_crl_hh}
    \Crl \bh_h \reff{eq_bh_II}= \suma \Crl \bha \reff{eq_bha_II}= \suma (\frh^\ver + \btha) \refff{eq_sigma_II}{eq_td_II}= \frh + \bdl_h.
\ee
This will be sufficient here, as in Theorem~\ref{thm_hp}, the minimization does not have a constrained curl.
It is also to be noted that the minimizations in~\eqref{eq_tha_II} avoid the additional elementwise orthogonality constraint of~\eqref{eq_tha}. Commuting could be achieved by further lowering the degree $p-1$ to $p-2$ and constructing an appropriate $\bdla$, but the above \cred{simplifications will be sufficient and actually} helpful for us in deriving $p$-robust bounds.

Definition~\ref{def_Phc_II} is well-posed, which can be easily verified as in Section~\ref{sec_Phc_WP}; namely, $\frh^\ver + \btha \in \RT_p(\Ta) \cap \Hdva$ with $\Dv (\frh^\ver + \btha )=0$. We now modify the reasoning of Step~2b from the proof of Lemma~\ref{lem_approx}, extending to the $\oma$-patch setting the $hp$-reasoning of Lemmas~\ref{lem_hp_d_osc} and~\ref{lem_hp_d_est} from Appendix~\ref{sec_constr} below.

\bl[$p$-robust patch estimate] \label{lem_p_rob_patch}
Let $\tv \in \HcN$ and a vertex $\ver \in \Vh$ be given. Let $\frh^\ver$ and $\bha$ be respectively given by Definitions~\ref{def_Phd_II} and~\ref{def_Phc_II}. \cred{If} $\GD = \emptyset$ and the patch subdomains $\oma$ are convex for all vertices $\ver \in \Vh$\cred{, set $\delta \eq 1$, otherwise set $\delta \eq 0$.} Then
\[
    \norm{\psia \biot_h - \bha}_\oma^2 \ls \sum_{\elm\in\Ta} \bigg[\min_{\tv_h \in \ND_{p-1}(\elm)}\norm{\tv - \tv_h}_\elm^2 + \bigg(\frac{h_\elm}{p^{\cred{\delta}}} \norm{\Crl\tv - \RToproj{p-1}(\Crl\tv)}_\elm\bigg)^2\bigg],
\]
where the hidden constant only depends on the shape-regularity parameter $\kappa_{\Th}$ of the mesh $\Th$.
\el

\bp Taking into account definition~\eqref{eq_bioth_II} and the constrained--unconstrained equivalence established in Lemma~\ref{lem_loc_const_eq} (employed with $p-1$ in place of $p$), we only need to show
\be \label{eq_p_rob_patch_proof}
    \norm{\psia \biot_h - \bha}_\oma \ls \norm{\tv - \biot_h}_\oma + \Bigg\{\sum_{\elm\in\Ta} \bigg(\frac{h_\elm}{p^{\cred{\delta}}} \norm{\Crl\tv - \RToproj{p-1}(\Crl\tv)}_\elm\bigg)^2\Bigg\}^\ft.
\ee
From the definition~\eqref{eq_bha_II}, as in~\eqref{eq_stab_patch},
\cite[Theorem~3.3]{Chaum_Voh_p_rob_3D_H_curl_23} gives
\ban
\norm{\psia \biot_h - \bha}_\oma & \ls \min_{\substack{
\bvf \in \Hca
\\
\Crl \bvf = \frh^\ver + \btha
}}
\norm{\psia \biot_h - \bvf}_\oma\\
& = \sup_{\substack{\bvf \in \Hcdual{\oma}\\ \norm{\Crl \bvf}_\oma = 1}}
    \big\{(\frh^\ver + \btha, \bvf)_\oma - (\psia \biot_h, \Crl \bvf)_\oma\big\},
\ean
where we have, for the equality, employed the primal--dual equivalence as in~\eqref{eq_stab_est}. Fix now
\be \label{eq_bvf}
    \bvf \in \Hcdual{\oma} \text{ with } \norm{\Crl \bvf}_\oma = 1.
\ee

Let $q \in \Hsa$ be such that
\be \label{eq_q_oma_2}
    (\Gr q, \Gr v)_\oma = (\cred{\bvf}, \Gr v)_\oma \qquad \forall v \in \Hsa
\ee
and define
\be \label{eq_bchi}
    \bchi \eq \bvf - \Gr q.
\ee
Recalling the notation~\cred{\eqref{eq_loc_spaces_int} or~\eqref{eq_loc_spaces_Dir}}, we have
\be \label{eq_bchi_reg_oma}
    \bchi \in \Hcdual{\oma} \cap \Hdva \text{ with } \Dv \bchi = 0
\ee
together with
\be \label{eq_bchi_prop}
    \Crl \bchi = \Crl \bvf.
\ee
\cred{Moreover, from~\eqref{eq_bchi_reg_oma}, \eqref{eq_bchi_prop}, and~\eqref{eq_bvf}, the Poincar\'e--Friedrichs--Weber inequality~\eqref{eq_weber_patch} gives
\be \label{eq_weber_bchi}
    \norm{\bchi}_\oma \reff{eq_weber_patch}\leq C_{\mathrm{PFW}} h_{\oma} \norm{\Crl \bchi}_{\oma} \reff{eq_bchi_prop}= C_{\mathrm{PFW}} h_{\oma} \norm{\Crl \bvf}_{\oma} \reff{eq_bvf}\ls h_{\oma}.
\ee
In this way, we will be gaining the $h_\elm$ factor in the second term of~\eqref{eq_p_rob_patch_proof}, but not the factor $h_\elm / p$.}

\cred{If $\GD = \emptyset$, we note that there is always zero normal trace prescribed on $\pt \oma$ in the space $\Hdva$, see~\eqref{eq_Hdva}, and, congruently, $\Hcdual{\oma} = \Hci{\oma}$, see~\eqref{eq_Hca_dual}. Thus, if additionally the patch subdomains $\oma$ are convex for all vertices $\ver \in \Vh$, we can employ~\cite[Theorem~3.9]{Gir_Rav_NS_86}. Thus, in this case,} there also holds $\bchi \in \tHoi{\oma}$ with
\be \label{eq_gr_chi_oma}
    \norm{\Gr \bchi}_\oma \leq \norm{\Crl \bchi}_{\oma} \reff{eq_bchi_prop}= \norm{\Crl \bvf}_{\oma} \reff{eq_bvf}=1.
\ee

Let
\be \label{eq_T}
    T \eq (\frh^\ver + \btha, \bvf)_\oma - (\psia \biot_h, \Crl \bvf)_\oma.
\ee
As in~\eqref{eq_vol}, since $\frh^\ver + \btha \in \RT_{p}(\Ta) \cap \Hdva$ is divergence-free, there exists $\ttw_h \in \ND_p(\Ta) \cap \Hca$ such that $\Crl \ttw_h = \frh^\ver + \btha$, and we can exchange $\bvf$ for \cred{$\bchi$ in~\eqref{eq_T}, \ie,
\ban
    T & = (\Crl \ttw_h, \bvf)_\oma - (\psia \biot_h, \Crl \bvf)_\oma = (\ttw_h, \Crl \bvf)_\oma - (\psia \biot_h, \Crl \bvf)_\oma \\
    & = (\ttw_h, \Crl \bchi)_\oma - (\psia \biot_h, \Crl \bchi)_\oma = (\frh^\ver + \btha, \bchi)_\oma - (\psia \biot_h, \Crl \bchi)_\oma.
\ean}
We will now work on this term.

Using the chain rule
\be \label{eq_chain}
    \Crl (\psia \bchi) = \Gr \psia \vp \bchi + \psia \Crl \bchi,
\ee
and since $(\cred{\td} \vp \cred{\te}) \scp \cred{\tf} = \cred{\te} \scp (\cred{\tf} \vp \cred{\td})$ and $\cred{\tf} \vp \cred{\td} = - \cred{\td} \vp \cred{\tf}$ for vectors $\cred{\td}, \cred{\te}, \cred{\tf} \in \RR^3$, we arrive at
\ban
    T = {} & (\frh^\ver + \btha, \bchi)_\oma + (\biot_h, \Gr \psia \vp \bchi)_\oma - (\biot_h, \Crl (\psia \bchi))_\oma \\
    = {} & (\frh^\ver + \btha - \Gr \psia \vp \biot_h, \bchi)_\oma - (\biot_h, \Crl (\psia \bchi))_\oma.
\ean
Moreover, the Green theorem gives
\[
    (\Crl \tv, (\psia \bchi))_\oma - (\tv, \Crl (\psia \bchi))_\oma = 0,
\]
so that altogether, also adding and subtracting $\psia \btau_h$,
\ban
    T = {} & (\frh^\ver + \btha - \Gr \psia \vp \biot_h - \psia \Crl \tv, \bchi)_\oma - (\biot_h - \tv, \Crl (\psia \bchi))_\oma \\
    = {} & (\frh^\ver - \psia \btau_h + \btha - \Gr \psia \vp \biot_h, \bchi)_\oma + (\btau_h - \Crl \tv, \psia \bchi)_\oma\\
    {} & - (\biot_h - \tv, \Crl (\psia \bchi))_\oma.
\ean
In the remainder of the proof, we estimate the above terms separately, in four steps.

{\em Step~1}. The chain rule~\eqref{eq_chain}, the Cauchy--Schwarz inequality, mesh shape-regularity, and~\eqref{eq_weber_bchi}, \cred{\eqref{eq_bchi_prop}, \eqref{eq_bvf}} give
\be \label{eq_psi_curl}
    \norm{\Crl (\psia \bchi))}_\oma = \norm{\Gr \psia \vp \bchi + \psia \Crl \bchi}_\oma \leq
    \norm{\Gr \psia}_{\infty,\oma} \norm{\bchi}_\oma + \norm{\Crl \bchi}_\oma \ls 1,
\ee
so that the Cauchy--Schwarz inequality yields
\[
    - (\biot_h - \tv, \Crl (\psia \bchi))_\oma \ls \norm{\biot_h - \tv}_\oma\cred{,}
\]
\cred{which is the first term in~\eqref{eq_p_rob_patch_proof}.}

{\em Step~2}. \cred{In general
\[
    (\btau_h - \Crl \tv, \psia \bchi)_\oma \leq \norm{\btau_h - \Crl \tv}_\oma \norm{\psia \bchi}_\oma \reff{eq_weber_bchi}\ls \norm{\btau_h - \Crl \tv}_\oma h_{\oma}
\]
and
\be \label{eq_th_est}
    \norm{\btau_h - \Crl \tv}_\oma^2 = \sum_{\elm\in\Ta} \norm{\btau_h - \Crl \tv}_\elm^2 \ls \sum_{\elm\in\Ta} \norm{\Crl\tv - \RToproj{p-1}(\Crl\tv)}_\elm^2,
\ee
relying on~\eqref{eq_tauh_II} and the $p$-robust constrained--unconstrained equivalence of~\cite[Lemma~A1]{Ern_Gud_Sme_Voh_loc_glob_div_22} (recall that $\Dv(\Crl \tv)=0$) (\cf~\eqref{eq_proj_est} below). Since $h_\oma \approx h_\elm$ in the patch, this gives the second term in~\eqref{eq_p_rob_patch_proof} with $\delta = 0$.}

\cred{If $\GD = \emptyset$ and the patch subdomains $\oma$ are convex for all vertices $\ver \in \Vh$, we rather} use~\eqref{eq_weber_bchi} and~\eqref{eq_gr_chi_oma}, which give, \cred{together with the chain rule,}
\[
    \norm{\Gr (\psia \bchi)}_\oma \cred{= \norm{\Gr \psia \bchi + \psia \Gr \bchi}_\oma \leq \norm{\Gr \psia}_{\infty, \oma} \norm{\bchi}_\oma + \norm{\Gr \bchi}_\oma} \ls 1.
\]
Thus, employing Lemma~\ref{lem_hp_d_est} for each $\elm \in \Ta$ and the Cauchy--Schwarz inequality, we infer
\ban
    (\btau_h - \Crl \tv, \psia \bchi)_\oma & = \sum_{\elm \in \Ta} (\btau_h - \Crl \tv, \psia \bchi)_\elm \\
    & \ls \sum_{\elm \in \Ta} \bigg( \frac{h_{\elm}}p\norm[\big]{\Crl\tv - \RToproj{p-1}(\Crl\tv)}_\elm \norm{\Gr(\psia \bchi)}_\elm \bigg) \\
    & \ls \Bigg\{\sum_{\elm\in\Ta} \bigg(\frac{h_\elm}{p} \norm{\Crl\tv - \RToproj{p-1}(\Crl\tv)}_\elm\bigg)^2\Bigg\}^\ft\cred{,}
\ean
\cred{which is the second term in~\eqref{eq_p_rob_patch_proof} with $\delta = 1$.}

{\em Step~3}. Relying on~\cite[Lemma~A.3]{Chaum_Voh_Maxwell_equil_23}, as in~\eqref{eq_tha_eff}, we have
\be \label{eq_tha_eff_II}
    \norm{\btha - \Gr \psia \vp \biot_h}_\oma \ls h_\oma^{-1} \norm{\tv - \biot_h}_\oma + h_\oma^{-1}\Bigg\{\sum_{\elm\in\Ta} \bigg(h_{\elm} \norm{\Crl\tv - \RToproj{p-1}(\Crl\tv)}_\elm \bigg)^2\Bigg\}^\ft;
\ee
note that this bound is indeed $p$-robust, as the above oscillation term is not of $hp$ type. \cred{Consequently,
\[
    (\btha - \Gr \psia \vp \biot_h, \bchi)_\oma \leq \norm{\btha - \Gr \psia \vp \biot_h}_\oma \norm{\bchi}_\oma,
\]
and we have~\eqref{eq_p_rob_patch_proof} with $\delta = 0$ from~\eqref{eq_weber_bchi}.}

\cred{If $\GD = \emptyset$ and the patch subdomains $\oma$ are convex for all vertices $\ver \in \Vh$, we again proceed more carefully.} The Euler--Lagrange conditions of~\eqref{eq_tha_II}, where it is important that no additional elementwise orthogonality constraints as in~\eqref{eq_tha} appear, give
\[
    (\btha - \Gr \psia \vp \biot_h, \bchi)_\oma = (\btha - \Gr \psia \vp \biot_h, \bchi - \bchi_h)_\oma
\]
for
\be \label{eq_bchi_oma}
    \bchi_h \eq \arg \min_{\substack{\ttw_h \in \RT_{p}(\Ta) \cap \Hdva\\
    \Dv \ttw_h = 0}}
    \norm{\bchi - \ttw_h}_\oma.
\ee
We now crucially rely on the $hp$-approximation estimate of~\cite[Theorem~3.6]{Ern_Gud_Sme_Voh_loc_glob_div_22} for $\Om = \oma$, recalling that $\bchi \in \Hdva$ has been prepared such that $\Dv \bchi = 0$ and $\bchi \in \tHoi{\oma}$, so that we can take $s=1$ in~\cite{Ern_Gud_Sme_Voh_loc_glob_div_22}. This gives
\be \label{eq_bchi_oma_est} \bs
    \norm{\bchi - \bchi_h}_\oma & \ls \Bigg\{\sum_{\elm\in\Ta} \bigg(\frac{h_\elm}{p+1} \big(\cred{h_\elm^{-1}} \norm{\bchi}_\elm^2 + \norm{\Gr \bchi}_\elm^2\big)^\ft\bigg)^2\Bigg\}^\ft \\
    & \ls \frac{h_\oma}{p} \big(\cred{h_\oma^{-1}}\norm{\bchi}_\oma^2 + \norm{\Gr \bchi}_\oma^2\big)^\ft \ls \frac{h_\oma}{p},
\es \ee
where we have employed~\cred{\eqref{eq_weber_bchi} and~\eqref{eq_gr_chi_oma}}.
Thus,
\[
    (\btha - \Gr \psia \vp \biot_h, \bchi)_\oma \ls \norm{\btha - \Gr \psia \vp \biot_h}_\oma \frac{h_\oma}{p}
\]
\cred{and we conclude~\eqref{eq_p_rob_patch_proof} with $\delta = 1$ using~\eqref{eq_tha_eff_II} (note that $p \geq 1$)}.

{\em Step~4}. \cred{We proceed as} above in Step~3. \cred{We employ~\eqref{eq_div_p_1} and either estimate
\[
    (\frh^\ver - \psia \btau_h , \bchi)_\oma \leq \norm{\frh^\ver - \psia \btau_h}_\oma \norm{\bchi}_\oma
\]
together with~\eqref{eq_weber_bchi}, or} we employ $\bchi_h$ from~\eqref{eq_bchi_oma} and~\eqref{eq_bchi_oma_est} \cred{together with} the Euler--Lagrange conditions of~\eqref{eq_sigma_a_II} and the Cauchy--Schwarz inequality, which give
\[
    (\frh^\ver - \psia \btau_h , \bchi)_\oma = (\frh^\ver - \psia \btau_h , \bchi - \bchi_h)_\oma \leq \norm{\frh^\ver - \psia \btau_h}_\oma \norm{\bchi - \bchi_h}_\oma.
\]
\cred{This concludes the proof.}
\ep

We are now ready to derive a $p$-robust variant of the inequality $\ls$ from~\eqref{eq_LG_unconstr} with a lowered polynomial degree on the right-hand side.

\begin{prop}[$p$-robust one-sided bound] \label{prop_p_rob_one_side}
Let $\tv \in \HcN$, a \cred{simplicial} mesh $\Th$ of $\Om$, and a polynomial degree $p \geq 1$ be fixed. \cred{If} $\GD = \emptyset$ and the patch subdomains $\oma$ are convex for all vertices $\ver \in \Vh$\cred{, set $\delta \eq 1$, otherwise set $\delta \eq 0$.}
Then
\be \label{eq_p_rob_all} \bs
{} & \min_{\tv_h \in \ND_p(\Th) \cap \HcN}\bigg[\norm{\tv - \tv_h}^2 + \sum_{\elm \in \Th}\bigg(\frac{h_{\elm}}{\cred{(}p+1\cred{)^\delta}}\norm{\Crl(\tv-\tv_h)}_\elm\bigg)^2\bigg] \\
\ls {} & \sum_{\elm\in\Th} \bigg[\min_{\tv_h \in \ND_{p-1}(\elm)}\norm{\tv - \tv_h}_\elm^2 + \bigg(\frac{h_{\elm}}{p^{\cred{\delta}}}\norm{\Crl\tv - \RToproj{p-1}(\Crl\tv)}_\elm\bigg)^2\bigg],
\es \ee
where the hidden constant only depends on the shape-regularity parameter $\kappa_{\Th}$ of the mesh $\Th$.
\end{prop}

\bp We bound the left-hand side of~\eqref{eq_p_rob_all} by
\[
\bigg[\norm{\tv - \bh_h}^2 + \sum_{\elm \in \Th}\bigg(\frac{h_{\elm}}{\cred{(}p+1\cred{)^\delta}}\norm{\Crl(\tv-\bh_h)}_\elm\bigg)^2\bigg],
\]
using $\bh_h$ from Definition~\ref{def_Phc_II}.
We now estimate the two above terms separately, in two steps.

{\em Step~1}. The triangle inequality gives, for $\biot_h$ given by~\eqref{eq_bioth_II},
\[
    \norm{\tv - \bh_h} \leq \norm{\tv - \biot_h} + \norm{\biot_h - \bh_h}.
\]
Thus, the constrained--unconstrained equivalence from Lemma~\ref{lem_loc_const_eq} (employed with $p-1$ in place of $p$) give\cred{s} the desired result for the first term above \cred{(with $\delta = 1$ in any case)}. As for the second one, definition~\eqref{eq_bh_II}, the partition of unity~\eqref{eq_PU}, and an overlap estimate imply
\[
    \norm{\biot_h - \bh_h}^2 = \sum_{\elm \in \Th} \norm[\Bigg]{\sum_{\ver \in \VK}(\psia \biot_h - \bha)}_\elm^2 \leq 4 \suma \norm{\psia \biot_h - \bha}_\oma^2,
\]
so that the desired bound follows from Lemma~\ref{lem_p_rob_patch} \cred{[--]}.

{\em Step~2}. From~\eqref{eq_crl_hh}, the partition of unity~\eqref{eq_PU}, \eqref{eq_sigma_II}, and adding and subtracting $\psia \Crl \biot_h$,
\ban
    \norm{\Crl(\tv-\bh_h)}_\elm & = \norm{\Crl\tv- \frh - \bdl_h}_\elm \\
    & = \norm[\Bigg]{\sum_{\ver \in \VK} \psia (\Crl(\tv - \biot_h)) + \sum_{\ver \in \VK} (\psia \Crl \biot_h - \frh^\ver) - \bdl_h}_\elm \\
    & \leq \sum_{\ver \in \VK} \norm{\psia (\Crl(\tv - \biot_h))}_\oma + \sum_{\ver \in \VK} \norm{\psia \Crl \biot_h - \frh^\ver}_\oma + \norm{\bdl_h}_\elm.
\ean
Recalling that the constraint in~\eqref{eq_bioth_II} gives $\Crl \biot_h = \btau_h$, we estimate the three terms separately.

For the first term, we estimate by~\eqref{eq_th_est}
\[
    \norm{\psia (\Crl(\tv - \biot_h))}_\oma \leq \norm{\btau_h - \Crl \tv}_\oma \reff{eq_th_est}\ls \Bigg\{\sum_{\elm\in\Ta} \norm{\Crl\tv - \RToproj{p-1}(\Crl\tv)}_\elm^2\Bigg\}^\ft.
\]
For the second one, we use~\eqref{eq_div_p_1} to see that
\[
    \norm{\psia \Crl \biot_h - \frh^\ver}_\oma \ls \Bigg\{\sum_{\elm\in\Ta} \norm{\Crl\tv - \RToproj{p-1}(\Crl\tv)}_\elm^2\Bigg\}^\ft.
\]
For the third one, we infer from~\eqref{eq_td_II} and~\eqref{eq_PU}
\[
    \norm{\bdl_h}_\elm = \norm[\Bigg]{\sum_{\ver \in \VK} (\btha - \Gr \psia \vp \biot_h)}_\elm \leq  \sum_{\ver \in \VK} \norm{\btha - \Gr \psia \vp \biot_h}_\oma.
\]
For the arising term, we can now use~\eqref{eq_tha_eff_II}. Thus, multiplying $\norm{\Crl(\tv-\bh_h)}_\elm$ by $h_{\elm}/\cred{(}p+1\cred{)^\delta}$, using that $h_{\elm}/\cred{(}p+1\cred{)^\delta} \leq h_{\elm}/p\cred{^\delta}$ as well as $1/\cred{(}p+1\cred{)^\delta} \leq 1$, and invoking the mesh shape-regularity, the proof is finished. \ep

Denote \label{page_proof_Thm_3_beg}
\be \label{eq_m}
    m^2 \eq \min_{\tv_h \in \ND_p(\Th) \cap \HcN}\bigg[\norm{\tv - \tv_h}^2 + \sum_{\elm \in \Th}\bigg(\frac{h_{\elm}}{p+1}\norm{\Crl(\tv-\tv_h)}_\elm\bigg)^2\bigg]
\ee
and
\be \label{eq_est_K}
    [v_{\elm,q,s,t}(\tv)]^2 \eq \bigg(\frac{h_\elm^{\min\{q,s\}}}{q^s} \norm{\tv}_{{\bm H}^s(\elm)}\bigg)^2 + \bigg(\frac{h_\elm}{q} \frac{h_\elm^{\min\{q,t\}}}{q^t} \norm{\Crl \tv}_{{\bm H}^t(\elm)} \bigg)^2.
\ee
From~\eqref{eq_LG_unconstr} in Theorem~\ref{thm_loc_glob} and the characterization~\eqref{eq_RToproj_eq}, we have
\be \label{eq_LG_unconstr_II}
m^2 \ls \sum_{\elm\in\Th} \bigg[\min_{\tv_h \in \ND_p(\elm)}\norm{\tv - \tv_h}_\elm^2 + \bigg(\frac{h_{\elm}}{p+1} \min_{\ttw_h \in \RT_p(\elm)} \norm{\Crl\tv - \ttw_h}_\elm\bigg)^2\bigg],
\ee
where the hidden constant can depend on $\kappa_{\Th}$ and the polynomial degree $p$.
From Proposition~\ref{prop_p_rob_one_side}, in turn,
\be \label{eq_p_rob_one_side_II}
m^2 \ls \sum_{\elm\in\Th} \bigg[\min_{\tv_h \in \ND_{p-1}(\elm)}\norm{\tv - \tv_h}_\elm^2 + \bigg(\frac{h_{\elm}}{p^{\cred{\delta}}} \min_{\ttw_h \in \RT_{p-1}(\elm)} \norm{\Crl\tv - \ttw_h}_\elm \bigg)^2\bigg],
\ee
where $p$ has to be greater or equal to $1$ but the hidden constant is independent of $p$. These are the two ingredients necessary for the proof of Theorem~\ref{thm_hp} in the spirit of the proof of Theorem~3.6 in~\cite[Section~5.2]{Ern_Gud_Sme_Voh_loc_glob_div_22}:

\bp[Proof of Theorem~\ref{thm_hp}, \cred{case~(ii) (convex patches and tangential boundary conditions}).]
\cred{Recall that the parameter $\delta$ from Proposition~\ref{prop_p_rob_one_side} is one here. Thus, employing} the notations~\eqref{eq_m}--\eqref{eq_est_K} and the claim~\eqref{eq_hp}, we need to show that
\be \label{eq_hp_proof}
    m^2 \ls \sum_{\elm \in \Th} [v_{\elm,p+1,s,t}(\tv)]^2,
\ee
where the hidden constant only depends on $\kappa_{\Th}$, $s$, and $t$. We distinguish two cases that we later combine together.

{\em Case~1}. Suppose $p \leq \cred{\max\{s,t\}}$. On each tetrahedron $\elm \in \Th$, we have from~\eqref{eq_N_K} and~\eqref{eq_RT_K} that $[\PP_p(\elm)]^3 \subset \ND_p(\elm)$ and $[\PP_p(\elm)]^3 \subset \RT_p(\elm)$. Thus the elementwise regularity assumptions~\eqref{eq_reg_v} and~\eqref{eq_reg_curl_v} together with~\cite[Lemma 4.1]{Bab_Sur_hp_FE_87} give
\[
    \min_{\tv_h \in \ND_p(\elm)}\norm{\tv - \tv_h}_\elm^2 + \bigg(\frac{h_{\elm}}{p+1} \min_{\ttw_h \in \RT_p(\elm)} \norm{\Crl\tv - \ttw_h}_\elm \bigg)^2 \ls [v_{\elm,p+1,s,t}(\tv)]^2,
\]
where the hidden constant only depends on $\kappa_{\Th}$, $s$, and $t$. Thus, from~\eqref{eq_LG_unconstr_II}, there exists a constant $C_{\kappa_{\Th}, p, s, t}$ only depending on $\kappa_{\Th}$, $p$, $s$, and $t$ such that
\[
    m^2 \leq C_{\kappa_{\Th}, p, s, t} \sum_{\elm \in \Th} [v_{\elm,p+1,s,t}(\tv)]^2.
\]
Defining $C_{\kappa_{\Th}, s, t}^\star \eq \max_{0 \leq p \leq \cred{\max\{s,t\}}} C_{\kappa_{\Th}, p, s, t}$, only depending on $\kappa_{\Th}$, $s$, and $t$, there holds, for all $0 \leq p \leq \cred{\max\{s,t\}}$,
\be \label{eq_est_1}
    m^2 \leq C_{\kappa_{\Th}, s, t}^\star \sum_{\elm \in \Th} [v_{\elm,p+1,s,t}(\tv)]^2.
\ee

{\em Case~2}. Suppose $p > \cred{\max\{s,t\}}$. Since $p$ is an integer, this implies that $p \geq 1$, so that we will be able to apply~\eqref{eq_p_rob_one_side_II}. As above, reducing $p+1$ to $p$ on the right-hand side,
\[
    \min_{\tv_h \in \ND_{p-1}(\elm)}\norm{\tv - \tv_h}_\elm^2 + \bigg(\frac{h_{\elm}}{p} \min_{\ttw_h \in \RT_{p-1}(\elm)} \norm{\Crl\tv - \ttw_h}_\elm \bigg)^2 \ls [v_{\elm,p,s,t}(\tv)]^2,
\]
where the hidden constant only depends on $\kappa_{\Th}$, $s$, and $t$.
Now, actually, $\min\{p,s\} = s = \min\{p+1,s\}$ and $p+1 \leq 2 p$, so that $1/p^s \leq 2^s/(p+1)^s$. Similarly, \cred{[--]} $\min\{p,t\} = t = \min\{p+1,t\}$ and $1/p^t \leq 2^t/(p+1)^t$. As a consequence, we can rise $p$ back to $p+1$ on the right-hand side,
\[
    \min_{\tv_h \in \ND_{p-1}(\elm)}\norm{\tv - \tv_h}_\elm^2 + \bigg(\frac{h_{\elm}}{p} \min_{\ttw_h \in \RT_{p-1}(\elm)} \norm{\Crl\tv - \ttw_h}_\elm \bigg)^2 \ls [v_{\elm,p+1,s,t}(\tv)]^2,
\]
where the constant beyond $\ls$ still only depends on $\kappa_{\Th}$, $s$, and $t$. Thus, using~\eqref{eq_p_rob_one_side_II} \cred{(recall that $\delta=1$ here)},
\be \label{eq_est_2}
    m^2 \leq C_{\kappa_{\Th}, s, t}^\sharp \sum_{\elm \in \Th} [v_{\elm,p+1,s,t}(\tv)]^2
\ee
for $C_{\kappa_{\Th}, s, t}^\sharp$ only depending on $\kappa_{\Th}$, $s$, and $t$ for all $p > \cred{\max\{s,t\}}$.

{\em Conclusion}. The claim\cred{~\eqref{eq_hp_proof} and consequently}~\eqref{eq_hp} follows from~\eqref{eq_est_1} and~\eqref{eq_est_2} for the constant $\max\{C_{\kappa_{\Th}, s, t}^\star, C_{\kappa_{\Th}, s, t}^\sharp\}$ only depending on $\kappa_{\Th}$, $s$, and $t$.
\ep

\cred{\bp[Proof of Theorem~\ref{thm_hp}, case~(i) (piecewise polynomial curl).] The proof is as above, replacing~\eqref{eq_est_K} by
\[
    [v_{\elm,q,s}(\tv)]^2 \eq \bigg(\frac{h_\elm^{\min\{q,s\}}}{q^s} \norm{\tv}_{{\bm H}^s(\elm)}\bigg)^2
\]
and considering the two cases $p \leq s$ and $p > s$. Since one assumes $\Crl \tv \in [\PP_{p-1}(\Th)]^3$, the second terms in~\eqref{eq_LG_unconstr_II} and~\eqref{eq_p_rob_one_side_II} vanish. \ep}

\cred{\bp[Proof of Theorem~\ref{thm_hp}, case~(iii) (general case).] The proof is as above, modulo two adjustments. First, we need to replace~\eqref{eq_est_K} by
\[
    [v_{\elm,q,s,t}(\tv)]^2 \eq \bigg(\frac{h_\elm^{\min\{q,s\}}}{q^s} \norm{\tv}_{{\bm H}^s(\elm)}\bigg)^2 + \bigg(h_\elm \frac{h_\elm^{\min\{q,t\}}}{q^t} \norm{\Crl \tv}_{{\bm H}^t(\elm)} \bigg)^2.
\]
Second, we use $h_{\elm}$ in place of $h_{\elm}/(p+1)$ in~\eqref{eq_m} and~\eqref{eq_LG_unconstr_II}; this is available from~\eqref{eq_LG_unconstr} in Theorem~\ref{thm_loc_glob}, noting that $1/(p+1) \leq 1$. Third, we need to employ~\eqref{eq_p_rob_one_side_II} from Proposition~\ref{prop_p_rob_one_side} with $\delta=0$. \label{page_proof_Thm_3_end}\ep}

\appendix

\section{$p$-robust equivalence of constrained and unconstrained best-appro\-xi\-ma\-ti\-on in \texorpdfstring{$\HC$}{H(curl)} on a tetrahedron} \label{sec_constr}

In this Appendix, we extend~\cite[Lemma~1]{Chaum_Voh_loc_glob_curl_21} to functions $\tv$ with nonpolynomial curl, in the spirit of~\cite[Lemma~A1]{Ern_Gud_Sme_Voh_loc_glob_div_22}. This is a consequence of the breakthrough result of Costabel and McIntosh in~\cite[Proposition~4.2]{Cost_McInt_Bog_Poinc_10}. Interestingly enough, the constant hidden in the inequality is here independent of the polynomial degree~$p$.

\bl[Equivalence of constrained and unconstrained best-approximation on a tetrahedron]\label{lem_loc_const_eq}
Let a polynomial degree $p \geq 0$, a tetrahedron $\elm$, and an arbitrary $\tv \in \Hci{\elm}$ be fixed. Let
\be \label{eq_tau_h}
    \btau_h \eq \arg \min_{\substack{\ttw_h \in \RT_p(\elm)
\\
\Dv \ttw_h = 0
}}
\norm{\Crl\tv - \ttw_h}_\elm.
\ee
Then
\be\label{eq_loc_const_eq} \bs
\min_{\tv_h \in \ND_p(\elm)}\norm{\tv - \tv_h}_\elm & \leq \min_{\substack{
\tv_h \in \ND_p(\elm)
\\
\Crl \tv_h = \btau_h
}}
\norm{\tv - \tv_h}_\elm \\
& \ls \Big(\min_{\tv_h \in \ND_p(\elm)}\norm{\tv - \tv_h}_\elm + \frac{h_{\elm}}{p+1}\norm[\big]{\Crl\tv - \RToproj{p}(\Crl\tv)}_\elm \Big),
\es \ee
where the hidden constant only depends on the shape-regularity $\kappa_\elm \eq h_\elm / \rho_\elm$ of $\elm$.
\el

\bp The first inequality is obvious, since the second minimization set has an additional curl constraint. In order to show the second one, denote respectively by
\be \label{eq_bioth_app}
    \biot_h \eq \arg \min_{\substack{
    \tv_h \in \ND_p(\elm)
    \\
    \Crl \tv_h = \btau_h
    }}
    \norm{\tv - \tv_h}_\elm
\ee
and
\be \label{eq_tbioth_app}
    \tilde \biot_h \eq \arg \min_{\tv_h \in \ND_p(\elm)}
    \norm{\tv - \tv_h}_\elm
\ee
the constrained and unconstrained minimizers. We then need to show
\be \label{eq_ineq}
    \norm{\tv - \biot_h}_\elm \ls \norm{\tv - \widetilde \biot_h}_\elm + \frac{h_{\elm}}{p+1}\norm[\big]{\Crl\tv - \RToproj{p}(\Crl\tv)}_\elm,
\ee
where $\ls$ means inequality up to a constant only depending on the shape-regularity $\kappa_\elm$.

Using~\eqref{eq_tau_h} and~\eqref{eq_curl_Np_map_loc}, $\btau_h - \Crl \tilde \biot_h \in \{\ttw_h \in \RT_p(\elm); \, \Dv \ttw_h = 0\}$.
Thus, we can use~\cite[Proposition~4.2]{Cost_McInt_Bog_Poinc_10}, \cf\ also the reformulation in~\cite[Theorem~2]{Chaum_Ern_Voh_curl_elm_20}, stipulating the existence of $\bvf_h \in \ND_p(\elm)$ with $\Crl \bvf_h = \btau_h - \Crl \tilde \biot_h$ such that
\be \label{eq_CM}
    \norm{\bvf_h}_\elm \ls \min_{\substack{\bvf \in \Hci{\elm}\\ \Crl \bvf = \btau_h - \Crl \tilde \biot_h}} \norm{\bvf}_\elm.
\ee
Shifting now the right-hand side of~\eqref{eq_CM} by $\tilde \biot_h$, we arrive at
\[
    \norm{\bvf_h}_\elm \ls \min_{\substack{\bvf \in \Hci{\elm}\\ \Crl \bvf = \btau_h}} \norm{\bvf - \tilde \biot_h}_\elm.
\]
A primal--dual equivalence as in,
\eg, \cite[Lemma~5.5]{Chaum_Ern_Voh_Maxw_22} implies (as in Section~\ref{sec_cont_sp_BC}, $\Hciz{\elm}$ is composed of those $\bvf \in \Hci{\elm}$ that verify $\bvf \vp \tn_\elm=0$ on  $\pt \elm$ in \cred{the weak} sense~\cred{\eqref{eq_tang_trace}})
\ban
\min_{\substack{\bvf \in \Hci{\elm}\\ \Crl \bvf = \btau_h}} \norm{\bvf - \tilde \biot_h}_\elm = {} & \sup_{\substack{\bvf \in \Hciz{\elm}\\ \norm{\Crl \bvf}_\elm = 1}}
\big\{(\btau_h, \bvf)_\elm - (\tilde \biot_h, \Crl \bvf)_\elm\big\} \\
\leq {} & \sup_{\substack{\bvf \in \Hciz{\elm}\\ \norm{\Crl \bvf}_\elm = 1}}
\big\{(\Crl \tv, \bvf)_\elm - (\tilde \biot_h, \Crl \bvf)_\elm\big\} \\
{} & + \sup_{\substack{\bvf \in \Hciz{\elm}\\ \norm{\Crl \bvf}_\elm = 1}}
(\btau_h - \Crl \tv, \bvf)_\elm \\
\ls {} & \min_{\substack{\bvf \in \Hci{\elm}\\ \Crl \bvf = \Crl \tv}} \norm{\bvf - \tilde \biot_h}_\elm + \frac{h_{\elm}}{p+1}\norm[\big]{\Crl\tv - \RToproj{p}(\Crl\tv)}_\elm,
\ean
where, to estimate the second term on the middle line, we have used the technical result of Lemma~\ref{lem_hp_d_osc} below.

Consequently, since $\tv \in \Hci{\elm}$ satisfies the curl constraint above,
\be \label{eq_CM_shift}
    \norm{\bvf_h}_\elm \ls \norm{\tv - \tilde \biot_h}_\elm + \frac{h_{\elm}}{p+1}\norm[\big]{\Crl\tv - \RToproj{p}(\Crl\tv)}_\elm.
\ee
Now note that $(\bvf_h + \tilde \biot_h) \in \ND_p(\elm)$ with $\Crl(\bvf_h + \tilde \biot_h) = \btau_h$. Thus, $\bvf_h + \tilde \biot_h$ belongs to the minimization set in~\eqref{eq_bioth_app}, and the minimization property~\eqref{eq_bioth_app} of $\biot_h$ implies $\norm{\tv - \biot_h}_\elm \leq \norm{\tv - (\bvf_h + \tilde \biot_h)}_\elm$. Thus, by virtue of the triangle inequality and using~\eqref{eq_CM_shift}, we altogether infer
\ban
    \norm{\tv - \biot_h}_\elm & \leq \norm{\tv - (\bvf_h + \tilde \biot_h)}_\elm \leq \norm{\tv - \tilde \biot_h}_\elm + \norm{\bvf_h}_\elm \\
    & \ls \norm{\tv - \tilde \biot_h}_\elm  + \frac{h_{\elm}}{p+1}\norm[\big]{\Crl\tv - \RToproj{p}(\Crl\tv)}_\elm,
\ean
\ie, \eqref{eq_ineq}, and the proof is finished.
\ep

\bl[$hp$ data oscillation]\label{lem_hp_d_osc} Let the assumptions of Lemma~\ref{lem_loc_const_eq} be verified. Then
\[
    \sup_{\substack{\bvf \in \Hciz{\elm}\\ \norm{\Crl \bvf}_\elm = 1}}
    (\btau_h - \Crl \tv, \bvf)_\elm \ls \frac{h_{\elm}}{p+1}\norm[\big]{\Crl\tv - \RToproj{p}(\Crl\tv)}_\elm,
\]
where the hidden constant only depends on the shape-regularity $\kappa_\elm \eq h_\elm / \rho_\elm$ of $\elm$. \el

\bp Fix $\bvf \in \Hciz{\elm}$ with $\norm{\Crl \bvf}_\elm = 1$. From~\eqref{eq_lift}, there exists $\bpsi \in \tHoi{\elm} \cap \Hciz{\elm}$ such that $\Crl \bpsi = \Crl \bvf$ and
\[
    \norm{\Gr \bpsi}_\elm \leq \norm{\Crl \bvf}_\elm = 1.
\]
Since $\btau_h \in \RT_p(\elm)$ with $\Dv \btau_h = 0$, there exists $\ttw_h \in \ND_p(\elm)$ such that $\Crl \ttw_h = \btau_h$. Thus, by the Green theorem,
\ban
    (\btau_h - \Crl \tv, \bvf)_\elm & = (\Crl (\ttw_h - \tv), \bvf)_\elm = (\ttw_h - \tv, \Crl \bvf)_\elm \\
    & = (\ttw_h - \tv, \Crl \bpsi)_\elm = (\btau_h - \Crl \tv, \bpsi)_\elm,
\ean
and we conclude by the following Lemma~\ref{lem_hp_d_est}.
\ep

\bl[$hp$ data estimate]\label{lem_hp_d_est} Let the assumptions of Lemma~\ref{lem_loc_const_eq} be verified. Let $\bpsi \in \tHoi{\elm}$. Then
\[
    \vert(\btau_h - \Crl \tv, \bpsi)_\elm\vert \ls \frac{h_{\elm}}{p+1}\norm[\big]{\Crl\tv - \RToproj{p}(\Crl\tv)}_\elm \norm{\Gr \bpsi}_\elm,
\]
where the hidden constant only depends on the shape-regularity $\kappa_\elm \eq h_\elm / \rho_\elm$ of $\elm$.
\el

\bp Let $q \in \Hooi{\elm}$ be such that
\be \label{eq_q}
    (\Gr q, \Gr v)_\elm = (\bpsi, \Gr v)_\elm \qquad \forall v \in \Hooi{\elm}.
\ee
Defining $\bchi \eq \bpsi - \Gr q$, we have
\be \label{eq_bchi_reg}
    \bchi \in \Hdvi{\elm} \text{ with } \Dv \bchi = 0.
\ee
Moreover, by the Green theorem, the right-hand side in~\eqref{eq_q} can be equivalently written as $(\bpsi, \Gr v)_\elm = - (\Dv \bpsi, v)_\elm$. Thus, since $\elm$ is convex, the elliptic regularity shift, see, \eg, \cite[\cred{Chapter~3}]{Gris_ell_nonsm_85}, gives $q \in \Hti{\elm}$ with
\[
    \norm{\Gr(\Gr q)}_\elm \ls \norm{\Dv \bpsi}_\elm \leq \norm{\Gr \bpsi}_\elm,
\]
where $\ls$ means inequality up to a constant only depending on the shape-regularity $\kappa_\elm$. Thus, in addition to~\eqref{eq_bchi_reg}, there also holds $\bchi \in \tHoi{\elm}$ with
\be \label{eq_gr_chi}
    \norm{\Gr \bchi}_\elm \leq \norm{\Gr \bpsi}_\elm + \norm{\Gr(\Gr q)}_\elm \ls \norm{\Gr \bpsi}_\elm.
\ee

Now, since from~\eqref{eq_tau_h} $\btau_h - \Crl \tv \in \Hdvi{\elm}$ with $\Dv (\btau_h - \Crl \tv) = 0$, there follows by the Green theorem
\[
    (\btau_h - \Crl \tv, \Gr q)_\elm = 0.
\]
Consequently,
\be \label{eq_II}
    (\btau_h - \Crl \tv, \bpsi)_\elm = (\btau_h - \Crl \tv, \bchi)_\elm.
\ee
Let respectively
\[
    \bchi_h \eq \arg \min_{\substack{\ttw_h \in \RT_p(\elm)\\
    \Dv \ttw_h = 0}}
    \norm{\bchi - \ttw_h}_\elm
\]
and
\[
    \tilde \bchi_h \eq \arg \min_{\ttw_h \in \RT_p(\elm)}
    \norm{\bchi - \ttw_h}_\elm
\]
be the constrained and unconstrained Raviart--Thomas approximations of $\bchi$.
The Euler--Lagrange conditions of~\eqref{eq_tau_h} allow us to subtract $\bchi_h$ (but not $\tilde \bchi_h$) in~\eqref{eq_II}, so that the Cauchy--Schwarz inequality and, crucially, the $p$-robust constrained--unconstrained $\HDV$ equivalence of~\cite[Lemma~A1]{Ern_Gud_Sme_Voh_loc_glob_div_22} (note that $\Dv \bchi = 0$) lead to
\be \label{eq_est} \bs
    (\btau_h - \Crl \tv, \bchi)_\elm
    & = (\btau_h - \Crl \tv, \bchi - \bchi_h)_\elm
    \leq \norm{\btau_h - \Crl \tv}_\elm \norm{\bchi - \bchi_h}_\elm \\
    & \ls \norm{\btau_h - \Crl \tv}_\elm \norm{\bchi - \tilde \bchi_h}_\elm.
\es \ee
Now, since $\RT_p(\elm)$ contains by~\eqref{eq_RT_K} polynomials up to degree $p$ in each component and since the minimization for $\tilde \bchi_h$ is unconstrained, the $hp$ approximation bound~\eqref{eq_Poinc_hp} gives
\be \label{eq_hp_Poinc}
    \norm{\bchi - \tilde \bchi_h}_\elm \ls \frac{h_{\elm}}{p+1} \norm{\Gr \bchi}_\elm.
\ee

Finally, in addition to~\eqref{eq_tau_h}, let
\[
    \tilde \btau_h \eq \arg \min_{\ttw_h \in \RT_p(\elm)} \norm{\Crl\tv - \ttw_h}_\elm,
\]
\ie, from~\eqref{eq_RToproj_eq}, $\tilde \btau_h = \RToproj{p}(\Crl\tv)$. It follows once again from the $p$-robust constrained--unconstrained $\HDV$ equivalence of~\cite[Lemma~A1]{Ern_Gud_Sme_Voh_loc_glob_div_22} (note that $\Dv (\Crl \tv) = 0$)  that
\be \label{eq_proj_est}
    \norm{\btau_h - \Crl \tv}_\elm \ls \norm{\tilde \btau_h - \Crl \tv}_\elm = \norm[\big]{\Crl\tv - \RToproj{p}(\Crl\tv)}_\elm.
\ee
Thus, the desired result is a combination of~\eqref{eq_II}, \eqref{eq_est}, \eqref{eq_hp_Poinc}, and~\eqref{eq_proj_est} together with~\eqref{eq_gr_chi}. \ep

\cred
{
\section{Broken regular decomposition in a patch} \label{sec_broken_reg}

Let $\Ta$ be a vertex patch for a mesh vertex $\ver \in \Vh$ as in Section~\ref{sec_mesh}. Recall the local spaces~\eqref{eq_loc_spaces_int} and~\eqref{eq_loc_spaces_Dir}, including the notation $\gD$ standing for the faces sharing the vertex $\ver$ and lying in $\overline \GD$. Then there holds:

\bl[Broken regular decomposition in a patch] \label{lem_broken_reg} Let $\bvf \in \Hcdual{\oma}$. Then there exists $\bpsi \in \Hcdual{\oma}$ with
\bse \label{eq_reg_dec} \be \label{eq_HTh}
    \bpsi\vert_\elm \in \tHoi{\elm} \qquad \forall \elm \in \Ta
\ee
%
such that
\be \label{eq_crl_eq}
    \Crl \bpsi = \Crl \bvf
\ee
and
\ba
    \norm{\bpsi}_\oma & \ls h_\oma \norm{\Crl \bvf}_\oma, \label{eq_est_L2}\\
    \Bigg\{\sum_{\elm \in \Ta} \norm{\Gr \bpsi}_\elm^2 \Bigg\}^\ft & \ls \norm{\Crl \bvf}_\oma, \label{eq_est_grad}
\ea \ese
where the hidden constants only depend on the shape-regularity parameter $\kappa_{\Th}$ of the mesh $\Th$.
\el

\bp
Our proof involves a reference patch. For a fixed shape-regularity parameter $\kappa_{\Th}$,
there is a maximal number of tetrahedra $\vert \Ta \vert$ in each vertex patch $\Ta$.
Therefore, there exists a finite number of ways these tetrahedra can be connected together
in the patch, and we may, for each value of $\kappa_{\Th}$, define a finite set
$\textup R$ of references patches of unit diameter in such a way that each vertex patch
$\Ta$ in the mesh is the image of an element $\widetilde {\mathcal{T}} \in \textup R$ through
a piecewise affine map that preserves connectivity.

Let $\widetilde {\mathcal{T}} \in \textup R$ be the reference patch associated with the given
$\Ta$, and denote by $\mathcal A: \widetilde \om \to \oma$ the corresponding piecewise affine map,
where $\widetilde \om$ is the open subdomain corresponding to $\widetilde {\mathcal{T}}$.
Then, for each $\elm \in \Ta$, we denote by
\begin{equation}
\label{eq_definition_piola}
\big (\mathcal A^{\rm c} \tv \big )\vert_\elm
\eq
\mathbb J^{\rm -T}_{\elm} \tv \circ \mathcal A^{-1},
\qquad
\big (\mathcal A^{\rm d} \tv\big )\vert_\elm
\eq
(\det \mathbb J_{\elm})^{-1} \mathbb J_{\elm} \tv \circ \mathcal A^{-1},
\end{equation}
the curl- and divergence-preserving Piola mappings associated with $\mathcal A$,
where $\mathbb J_{\elm}$ is the (constant) Jacobian matrix of $\mathcal A$ on $\elm$,
see e.g. \cite[Section 7.2]{Ern_Guermond_FEs_I_21}. It is standard that
$\mathcal A^{\rm c}$ maps $\Hcdual{\widetilde \om}$ into $\Hcdual{\oma}$,
where $\Hcdual{\widetilde \om}$ embeds essential boundary conditions on
$\mathcal A^{-1}(\gD)$. $\mathcal A^{\rm c}$ and $\mathcal A^{\rm d}$ are
bijective, and we denote by $(\mathcal A^{\rm c})^{-1}$ and $(\mathcal A^{\rm d})^{-1}$
their inverses. Finally, we have the commuting property
$\Crl (\mathcal A^{\rm c} \cdot) = \mathcal A^{\rm d}(\Crl \cdot)$.

Consider $\bvf \in \Hcdual{\oma}$ and let
$\widetilde \bvf \eq (\mathcal A^{\rm c})^{-1}(\bvf) \in \Hcdual{\widetilde \om}$.
We can employ~\eqref{eq_lift} (where we take $\ttw$ in $\tHoDi{\om}$) on the reference
patch $\widetilde \om$, where $\Clift{\widetilde \om}$ is now a generic constant only
depending on $\kappa_{\Th}$ due to the above discussion. Therefore, there exists
$\widetilde \bpsi \in \tHoi{\widetilde \om}$, with $\widetilde \bpsi = 0$ on
$\mathcal A^{-1}(\gD)$ if $\gD$ includes at least one face and of componentwise
mean value zero on $\widetilde \om$ otherwise such that
\begin{equation*}
\Crl \widetilde \bpsi = \Crl \widetilde \bvf \quad \text{ and } \quad
\quad
\norm{\Gr \widetilde \bpsi}_{\widetilde \om}
\leq
\Clift{\widetilde \om}
\norm{\Crl \widetilde \bvf}_{\widetilde \om}.
\end{equation*}
We also note that we have
\begin{equation}
\label{eq_poincare_reference}
\norm{\widetilde \bpsi}_{\widetilde \om}
\ls
\norm{\Gr \widetilde \bpsi}_{\widetilde \om}
\ls
\norm{\Crl \widetilde \bpsi}_{\widetilde \om}
=
\norm{\Crl \widetilde \bvf}_{\widetilde \om}
\end{equation}
by the Poincar\'e--Friedrichs inequality~\eqref{eq_Poinc_Fried}, since $\widetilde \om$ is of unit diameter.

We then introduce $\bpsi \eq \mathcal A^{\rm c}(\widetilde \bpsi) \in \Hcdual{\oma}$
and claim that it satisfies all the requirements in~\eqref{eq_reg_dec}.
Indeed, elementwise, the Piola mapping $\mathcal A^{\rm c}$ amounts to an affine
change of coordinates and a multiplication by a constant matrix. Therefore, it is
clear that it preserves smoothness elementwise and~\eqref{eq_HTh} follows from
$\widetilde \bpsi \in \tHoi{\widetilde \om}$. Moreover, \eqref{eq_crl_eq} holds true since
\begin{align*}
\Crl \bpsi
& =
\Crl \mathcal A^{\rm c}(\widetilde \bpsi)
=
\mathcal A^{\rm d}(\Crl \widetilde \bpsi)
=
\mathcal A^{\rm d}(\Crl \widetilde \bvf)
=
\mathcal A^{\rm d}(\Crl ((\mathcal A^{\rm c})^{-1}(\bvf)))\\
& =
\mathcal A^{\rm d}((\mathcal A^{\rm d})^{-1}(\Crl \bvf))
=
\Crl \bvf.
\end{align*}

The estimates~\eqref{eq_est_L2} and~\eqref{eq_est_grad} then follow from
\eqref{eq_poincare_reference} and the usual scaling arguments due to the
structure~\eqref{eq_definition_piola} of the Piola mappings.
Indeed, for each $\elm \in \Ta$, we have, on the one hand,
\ban
\norm{\bpsi}_{\elm}^2
& =
\norm{\mathcal A^{\rm c}(\widetilde \bpsi)}_{\elm}^2
=
\norm{\mathbb J_\elm^{-T}\widetilde \bpsi \circ \mathcal A^{-1}}_{\elm}^2
\\
& \lesssim
h_\elm^{-2} \norm{\widetilde \bpsi \circ \mathcal A^{-1}}_\elm^2
\lesssim
h_\elm^{-2} \vert \elm\vert \norm{\widetilde \bpsi}_{\widetilde \elm}^2
\lesssim
h_\elm \norm{\widetilde \bpsi}_{\widetilde \elm}^2,
\ean
where $\widetilde \elm = \mathcal A^{-1}(\elm)$. On the other hand, it is true that
\be \label{tmp_estimate_curl_tbpsi} \bs
\norm{\Crl \widetilde \bpsi}_{\widetilde \elm}^2
& =
\norm{(\mathcal A^{\rm d})^{-1}(\Crl \bpsi)}_{\widetilde \elm}^2
=
\norm{(\det \mathbb J_{\elm}) \mathbb J_\elm^{-1} (\Crl \bpsi) \circ \mathcal A}_{\widetilde \elm}^2
\\
& \lesssim
(\det \mathbb J_\elm)^2 h_\elm^{-2}
\norm{(\Crl \bpsi) \circ \mathcal A}_{\widetilde \elm}^2
\lesssim
\frac{(\det \mathbb J_\elm)^2 h_\elm^{-2}}{\vert \elm \vert}
\norm{\Crl \bpsi}_{\elm}^2\\
& \lesssim
h_\elm \norm{\Crl \bpsi}_\elm^2,
\es \ee
so that~\eqref{eq_est_L2} follows from~\eqref{eq_poincare_reference}.

Finally, we establish~\eqref{eq_est_grad} by directly differentiating in
\eqref{eq_definition_piola}. Specifically, we have
\begin{equation*}
\partial_\ell \bpsi
=
\mathbb J_\elm^{-T} \partial_\ell(\widetilde \bpsi \circ \mathcal A^{-1})
=
\mathbb J_\elm^{-T} \partial_\ell (\mathcal A^{-1}_j) (\partial_j\widetilde \bpsi \circ \mathcal A^{-1}),
\end{equation*}
so that
\begin{equation*}
\vert\partial_\ell \bpsi\vert
\leq
h_\elm^{-2} \vert\Gr \widetilde \bpsi \circ \mathcal A^{-1}\vert
\end{equation*}
and
\begin{equation*}
\norm{\partial_\ell \bpsi}_\elm^2
\lesssim
h^{-4}_\elm \vert \elm \vert \norm{\Gr \widetilde \bpsi}_{\widetilde \elm}^2.
\end{equation*}
We therefore arrive at
\begin{equation*}
\sum_{\elm \in \Ta} \norm{\Gr \bpsi}_{\elm}^2
\lesssim
h_{\oma}^{-1}\norm{\Gr \widetilde \bpsi}_{\widetilde \omega}^2,
\end{equation*}
and the conclusion follows from~\eqref{eq_poincare_reference} and
\eqref{tmp_estimate_curl_tbpsi}.
\ep
}

\subsection*{Acknowledgment}

This project has received funding from the European Research Council (ERC) under the European Union's Horizon 2020 research and innovation program (grant agreement No 647134 GATIPOR).

\ifNM\else\bibliographystyle{acm_mod}\fi
\bibliography{biblio}

\end{document}